\documentclass[12pt]{article}  
\def\sq{\hbox {\rlap{$\sqcap$}$\sqcup$}}
\def\sq{\hbox {\rlap{$\sqcap$}$\sqcup$}}
\def\R{ {\rm R \kern -.31cm I \kern .15cm}}
\def\C{ {\rm C \kern -.15cm \vrule width.5pt \kern .12cm}}
\def\Z{ {\rm Z \kern -.27cm \angle \kern .02cm}}
\def\N{ {\rm N \kern -.26cm \vrule width.4pt \kern .10cm}}
\def\1{{\rm 1\mskip-4.5mu l} }
\def\lsim{\raise0.3ex\hbox{$<$\kern-0.75em\raise-1.1ex\hbox{$\sim$}}}
\def\gsim{\raise0.3ex\hbox{$>$\kern-0.75em\raise-1.1ex\hbox{$\sim$}}}
\def\noi{\noindent}

\def\beq{\begin{equation}}   \def\eeq{\end{equation}}
\def\bea{\begin{eqnarray}}  \def\eea{\end{eqnarray}}
\def\nn{\nonumber}
\def\noi{\noindent}
\def\beeq{\begin{eqnarray}} \def\eeeq{\end{eqnarray}}
\newcommand\mysection{\setcounter{equation}{0}\section}

\newcounter{hran}

\begin{document} 
\centerline{\Large\bf Long Range Scattering for the  Modified} 
 \vskip 3 truemm \centerline{\Large\bf  Schr\"odinger Map in two space dimensions\footnote{Work supported in part by the EC contract RITA-CT-2004-505493}} \vskip 0.8 truecm

\centerline{\bf J. Ginibre}
\centerline{Laboratoire de Physique Th\'eorique\footnote{Unit\'e Mixte de
Recherche (CNRS) UMR 8627}}  \centerline{Universit\'e de Paris XI, B\^atiment
210, F-91405 ORSAY Cedex, France}
\vskip 3 truemm

\centerline{\bf G. Velo}
\centerline{Dipartimento di Fisica, Universit\`a di Bologna}  \centerline{and INFN, Sezione di
Bologna, Italy}

\vskip 1 truecm

\begin{abstract}
We study the asymptotic behaviour in time of solutions and the theory of scattering for the modified Schr\"odinger map in two space dimensions. We solve the Cauchy problem with large finite initial time, up to infinity in time, and we determine the asymptotic behaviour in time of the solutions thereby obtained. As a by product, we obtain global existence for small data in $H^k \cap FH^k$ with $k > 1$. We also solve the Cauchy problem with infinite initial time, namely we construct solutions defined in a neighborhood of infinity in time, with prescribed asymptotic behaviour of the previous type.
\end{abstract}

\vskip 1.5 truecm
\noi MS Classification :    Primary 35P25. Secondary 35B40, 35Q60.\par \vskip 2 truemm

\noi Key words : Long range scattering, modified Schr\"odinger map. \par 
\vskip 1 truecm

\noindent LPT Orsay 06-59\par
\noindent October 2006\par \vskip 3 truemm

\newpage
\pagestyle{plain}
\baselineskip 18pt

\mysection{Introduction}
\hspace*{\parindent} This paper is devoted to the study of the asymptotic behaviour in time of solutions and to the theory of scattering for the modified Schr\"odinger map (MSM) system in space dimension 2. In general space dimension $n$, that system takes the form
\beq
\label{1.1e}
i \partial_t u = - (1/2) \Delta_A u + g(u) u\ .
\eeq

\noi Here $u$ is a $\C^n$ vector valued function defined in space time ${I\hskip-1truemm R}^{n+1}$, $\Delta_A = \nabla_A^2 = (\nabla - i A)^2$ is the covariant Laplacian associated with the vector potential $A$ defined by
\beq
\label{1.2e}
A_j = 4 \Delta^{-1}\  \partial_k \ {\rm Im}\ \overline{u}_k \ u_j \ ,
\eeq

\noi $g(u)$ is the hermitian matrix defined by 
\bea
\label{1.3e}
&&g_{jk}(u) = - 2i\ {\rm Im} \ \overline{u}_k \ u_j - A_0\ \delta_{jk}\ ,\\
&&A_0 = 2 \Delta^{-1} \ \partial_j \ \partial_k {\rm Re}\ \overline{u}_k\ u_j - |u|^2
\label{1.4e}
\eea

\noi and summation over repeated indices is understood. The normalization and sign conventions in (\ref{1.1e})-(\ref{1.4e}) differ from those currently used by a few signs and factors of 2 in order to allow for an easier comparison with the Maxwell-Schr\"odinger (MS) system. The MSM system is formally derived from the more primitive  system \cite{16r}
\beq
\label{1.5e}
i \partial_t z = - (1/2) \nabla_j \partial_j z
\eeq

\noi where $z$ is a complex function defined in space time ${I\hskip-1truemm R}^{n+1}$ and
\beq
\label{1.6e}
\nabla_j = \partial_j - 2 \left ( 1 + |z|^2\right )^{-1}\ \overline{z} \left ( \partial_j z\right )\ .
\eeq

The system (\ref{1.5e}) itself is obtained through a stereographic projection from a more geometrically defined Schr\"odinger map (SM) system where the unknown function takes values in the unit sphere $S^2$. The latter system appears as the Landau-Lifschitz model of a ferromagnet \cite{14r}. When deriving the MSM system from the SM one, in addition to  (\ref{1.1e})-(\ref{1.4e}), one obtains a constraint satisfied by $u$. That constraint is easily seen to be formally preserved by the evolution (\ref{1.1e}). The MSM system can be studied with or without that constraint. In the present paper we consider it without the constraint, which makes it more general and significantly different as regards scattering (see below). It is an important question to make the correspondence between the SM and MSM systems rigorous, in order to transfer results from one system to the other. The equivalence of the SM system to the MSM one (with the constraint) has recently been proven under mild regularity assumptions \cite{NEW19r}.\par

A large amount of work has been devoted to the Cauchy problem both for the system (\ref{1.5e}) and for the MSM system (\ref{1.1e}) in various space dimensions. We refer to \cite{1r}-\cite{NEW5r}, \cite{7r}-\cite{10r}, \cite{12r}, \cite{NEW15r}, \cite{15r}-\cite{17r}, \cite{NEW24r} and the literature therein quoted. The best available result so far for the MSM system without the constraint in space dimension 2 is local wellposedness in $H^s$ for $s > 3/4$ \cite{10r}.\par

In the present paper we shall study the asymptotic behaviour in time of solutions and the theory of scattering for the MSM system without the constraint in space dimension 2, where it is borderline long range (see below). \par

Here we regard scattering theory as a method to classify the possible solutions of (\ref{1.1e}) by their asymptotic behaviour. That point of view leads to the following two problems.\\

\noi \underline{\bf Problem 1.} One gives oneself a set ${\cal U}_a$ of presumed asymptotic behaviours $u_a$ for the system (\ref{1.1e}), parametrized by some data $u_+$. For each $u_a \in {\cal U}_a$, one tries to construct a solution of the system (\ref{1.1e}) such that $u(t) - u_a(t)$ tends to zero as $t \to + \infty$ in a suitable sense, more precisely in suitable norms. The same problem can be considered for $t \to - \infty$. We restrict our attention to the case of $t \to + \infty$. The previous problem decomposes into two steps. The first step is to construct the solution $u$ in a neighborhood of $t= + \infty$, namely in an interval $[T, \infty )$ for $T$ sufficiently large. This is the local Cauchy problem at infinity in time. The second and rather independent step consists in extending the solution to all times and reduces therefore to the global Cauchy problem at finite times. In this paper we consider only the first step and leave aside the second one.\par

If the previous problem can be solved for any $u_a \in {\cal U}_a$, the map $u_a \to u$ thereby defined is essentially the wave operator $\Omega_+$ for positive time associated with ${\cal U}_a$.\\

 \noi \underline{\bf Problem 2.} This is the converse to Problem 1. Given a generic solution $u$ of the system (\ref{1.1e}), one tries to find an asymptotic motion $u_a \in {\cal U}_a$ such that $u(t) - u_a(t)$ tends to zero as $t \to \infty$ in a suitable sense (in suitable norms). If that problem and the same one for $t \to - \infty$ can be solved for all $u$ (in a suitable functional framework), one says that asymptotic completeness holds with respect to the set ${\cal U}_a$. That property requires in particular that all possible asymptotic behaviours of the solutions of (\ref{1.1e}) have been identified and included in ${\cal U}_a$, and is completely out of reach in the present case. In this paper, we restrict ourselves to the construction of a set of solutions of the Cauchy problem with finite initial time, defined up to infinity in time and behaving asymptotically as functions in the set ${\cal U}_a$ for which we can solve the first problem. That set includes small global solutions of (\ref{1.1e}). \\
 
 In the present case, a natural candidate for ${\cal U}_a$ is the set of solutions of the free Schr\"odinger equation, namely
\beq
\label{1.7e}
u_a(t) = U(t)\ u_+ = \exp \left ( i (t/2)\Delta \right ) u_+\  .
\eeq

\noi Cases where such a choice is adequate are referred to as short range cases. This requires the nonlinear interaction to decrease sufficiently fast at infinity in space and/or time. This occurs for the MSM system in space dimension $n \geq 3$ and, if the constraint is included, also in space dimension n = 2. This also occurs for the SM system in space dimension $n \geq 2$ \cite{NEW24r}. The long range case is the complementary one where that set is inadequate and has to be replaced by a set of modified asymptotic behaviours. This occurs for the MSM system without the constraint in space dimension 2, the case which we treat in this paper. The modification includes the introduction of a phase in the asymptotic Schr\"odinger function. In that respect, the MSM system without the constraint in space dimension 2 is borderline long range and similar to the Maxwell-Schr\"odinger (MS) and Wave-Schr\"odinger (WS) systems in dimension 3, and to the Hartree equation with $|x|^{-1}$ potential in dimension $n \geq 2$. \par

The MSM system is also similar to the MS system in the sense that it consists of a Schr\"odinger equation in a magnetic field, with vector potential $A$. In contrast with the MS system however, the magnetic field does not propagate, and the vector potential is defined locally in time in terms of the Schr\"odinger function. This makes the problem simpler and makes it possible to apply the methods previously used for the MS system in dimension 3 to the MSM system in dimension 2. \par

The theory of scattering for the MS system in dimension 3 has been studied by several authors \cite{4r} \cite{5r} \cite{6r} \cite{18r} \cite{20r}, following work on the WS system and on the Hartree equation. We refer to \cite{6r} for additional information and references on that matter. In this paper we study the MSM system in dimension 2 by the methods used in \cite{4r}, which seem to be the most readily applicable to that system. The main results are as follows. We first solve the Cauchy problem with large finite initial time, up to infinity in time and we determine the asymptotic behaviour in time of the solutions thereby obtained. This represents our contribution to Problem 2 mentioned above, and allows us to identify a set ${\cal U}_a$ of possible asymptotic behaviours. As a by product we obtain global existence for small data. We then solve the Cauchy problem with infinite initial time, namely we construct solutions defined in a  neighborhood of infinity in time, with prescribed asymptotic behaviour of the previous type. The method consists in expressing the Schr\"odinger function $u$ in terms of a complex amplitude $v$ and a real phase $\varphi$, replacing the original system (\ref{1.1e}) by an auxiliary system for the pair $(v, \varphi )$, treating the corresponding problems for the latter system, and reconstructing the solutions $u$ of the original system from the solutions $(v, \varphi )$ of the auxiliary one. The detailed construction is too complicated to allow for a more precise description at this stage and will be described in heuristic terms in Section 2 below. At the end of that section, we shall also give a simplified version of the results as Propositions 2.1 and 2.2. We conclude this introduction by giving a brief outline of the contents of this paper. A more detailed description will be given at the end of Section 2. In Section 3, we collect some notation and preliminary estimates. In Section 4 we study the Cauchy problem at finite initial time both for the original  system (\ref{1.1e}) and for the auxiliary system. In Section 5, we study the Cauchy problem with infinite initial time for the original system (\ref{1.1e}) and the corresponding problem for the auxiliary system.

\mysection{Heuristics and formal computations}
\hspace*{\parindent} In this section, we perform in a formal way the algebraic computations needed to study the Cauchy problem for the MSM system (\ref{1.1e}) in a neighborhood of infinity in time, both for finite and infinite initial time, and we sketch the method used to solve that problem. The system (\ref{1.1e}) in this form is not well suited for that purpose and we perform a number of transformations leading to an auxiliary system for which that problem can be handled. The unitary group 
\beq
\label{2.1e}
U(t) = \exp (i(t/2)\Delta )
\eeq

\noi which solves the free Schr\"odinger equation can be written as
\beq
\label{2.2e}
U(t) = M(t)\ D(t)\ F\ M(t)
\eeq

\noi where $M(t)$ is the operator of multiplication by the function 
\beq
\label{2.3e}
M(t) = \exp \left ( ix^2/2t \right ) \ ,
\eeq

\noi $F$ is the Fourier transform and $D(t)$ in the dilation operator defined by
\beq
\label{2.4e}
D(t) = (it)^{-1}\ D_0(t)\qquad , \quad \left ( D_0(t) f\right ) (x) = f(x/t)\ .
\eeq

\noi We first change variables from $u$ to its pseudoconformal inverse $u_c$ defined by 
\beq
\label{2.5e}
u(t) = M(t)\ D(t) \ \overline{u_c(1/t)}
\eeq

\noi or equivalently
\beq
\label{2.6e}
\widetilde{u}(t) =  \overline{F\widetilde{u}_c(1/t)}
\eeq

\noi where for any function $f$ of space time, we define 
\beq
\label{2.7e}
\widetilde{f}(t, \cdot) = U(-t) f(t, \cdot)\ .
\eeq

\noi Correspondingly we define $B$ by 
\beq
\label{2.8e}
A(t) = - t^{-1} D_0(t)\ B(1/t).
\eeq

\noi Substituting (\ref{2.5e}) (\ref{2.8e}) into (\ref{1.1e}) yields the following equation for $u_c$~:
\beq
\label{2.9e}
i \partial_t u_c = - (1/2) \Delta_B u_c + \left ( {\check B} (u_c) + g(u_c)\right )u_c
\eeq

\noi where $B = B(u_c)$ and we have defined
\beq
\label{2.10e}
B(v) = B(v, v)\quad , \quad {\check B}(v) = {\check B}(v, v)\quad , \quad g(v) = g(v,v)\ ,
\eeq
\beq
\label{2.11e}
B_j (v_1, v_2) = 2 \Delta^{-1} \partial_k \ {\rm Im} \left ( \overline{v}_{1k} v_{2j} + \overline{v}_{2k} v_{1j}\right ) \ ,
\eeq
\beq
\label{2.12e}
g_{jk}(v_1,v_2) = - i \ {\rm Im} \left ( \overline{v}_{1k} v_{2j} + \overline{v}_{2k} v_{1j}\right ) - B_0 (v_1, v_2) \delta_{jk}\ ,
\eeq
\beq
\label{2.13e}
B_0(v_1,v_2) = 2 \Delta^{-1} \partial_j \partial_k\ {\rm Re}  \left ( \overline{v}_{1k} v_{2j}\right ) - {\rm Re}  \left ( \overline{v}_{1} \cdot v_2 \right ) \ ,
\eeq
\beq
\label{2.14e}
{\check B}(v_1,v_2) = 2 t^{-1} x_j \Delta^{-1}  \partial_k\ {\rm Im}  \left ( \overline{v}_{1k} v_{2j} + \overline{v}_{2k} v_{1j} \right ) 
\eeq

\noi and more generally, for any ${I\hskip-1truemm R}^{2}$ vector valued function of space time 
\beq
\label{2.15e}
{\check f}(x,t) = t^{-1} x \cdot f(x, t)\ .
\eeq

\noi We next parametrize $u_c$ in terms of a complex amplitude $v$ and a real phase $\varphi$ by
\beq
\label{2.16e}
u_c = v \exp (- i \varphi ) \ .
\eeq

\noi Note that $u_c$ and $v$ are ($\C^2$) vector valued and that the phase is the same for all components, so that in particular
\beq
\label{2.17e}
B(u_c) = B(v) \quad , \quad {\check B}(u_c) =  {\check B}(v)\quad , \quad g(u_c) = g(v)\ .
\eeq

\noi Substituting (\ref{2.16e}) into (\ref{2.9e}) yields the equation 
\beq
\label{2.18e}
i \partial_t v = - (1/2) \Delta_K v + \left (  {\check B}(v) - \partial_t \varphi + g(v) \right ) v
\eeq

\noi where
\beq
\label{2.19e}
K = s + B\quad , \qquad s = \nabla \varphi 
\eeq

\noi and in the same way as before
$$\Delta_K = \nabla_K^2 \qquad ,\quad \nabla_K = \nabla - i K \ .$$

\noi We have now only one equation for two functions $(v, \varphi )$. We then arbitrarily impose a second equation, namely an equation for the phase $\varphi$, thereby splitting (\ref{2.18e}) into a system of two equations, the other one of which is an equation for $v$. There is a large amount of freedom in the choice of the equation for the phase. The role of the phase is to cancel the long range term $ {\check B}(v)$ in (\ref{2.18e}). However since that term has a relatively low regularity, it is convenient to split it into a short range and a long range part. Let $\chi \in {\cal C}^{\infty} ({I\hskip-1truemm R}^{2}, {I\hskip-1truemm R})$, $0 \leq \chi \leq 1$, $\chi (\xi ) = 1$ for $|\xi| \leq 1$, $\chi (\xi) = 0$ for $|\xi| \geq 2$. We define 
\beq
\label{2.20e}
{\check B}_L = F^*\chi (\ \cdot \ t^{1/2}) F {\check B} \qquad , \quad {\check B}_S = F^*\left ( (1 - \chi (\ \cdot \ t^{-1/2})\right )  F {\check B} \ .
\eeq

\noi As the equation for $\varphi$, we take
\beq
\label{2.21e}
\partial_t \varphi = {\check B}_L (v)
\eeq

\noi so that the equation for $v$ becomes
\beq
\label{2.22e}
i \partial_t v = Hv
\eeq

\noi with 
\beq
\label{2.23e}
H = - (1/2) \Delta_K + {\check B}_S(v) + g(v)\ .
\eeq

\noi The system (\ref{2.21e}) (\ref{2.22e}) is the final form of the auxiliary system that replaces the original system (\ref{1.1e}). For technical reasons, it will be useful to consider also the partly linearized system for a new variable $v'$  
\beq
\label{2.24e}
i \partial_t v' = H v'
\eeq

\noi where $H$ is still associated with $(v, \varphi )$ according to (\ref{2.23e}). The Cauchy problem at finite initial time $t_0 \in [1, \infty )$ for the original system (\ref{1.1e}) is now replaced by the Cauchy problem for the auxiliary system (\ref{2.21e}) (\ref{2.22e}) at finite initial time $\tau_0 = t_0^{-1} \in (0, 1]$. We shall solve that problem in two steps. We shall first solve the linearized system (\ref{2.24e}) for $v'$ with given $v$, thereby defining a map $\Gamma : v\to v'$. We shall then prove that the map $\Gamma$ has a fixed point by a contraction method. With the solution of the system (\ref{2.21e}) (\ref{2.22e}) available, the original system (\ref{1.1e}) can be solved by substituting the solution $(v , \varphi )$ into the formulas (\ref{2.6e}) (\ref{2.16e}).

In a similar way, the Cauchy problem at infinite initial time for the original system (\ref{1.1e}) is replaced by the Cauchy problem at $t=0$ for the auxiliary system (\ref{2.21e}) (\ref{2.22e}). Since that system is singular at $t=0$, that problem cannot be treated directly and we follow instead an indirect procedure. We choose a set of asymptotic functions $(v_a, \varphi_a)$ which are expected to be suitable asymptotic forms of $(v, \varphi )$ at $t = 0$ and we try to construct solutions of the auxiliary system (\ref{2.21e}) (\ref{2.22e}) that are asymptotic to $(v_a, \varphi_a)$ as $t \to 0$. The set $(v_a, \varphi_a)$ will be taken of the following form. For a given $v_a$, which needs not be specified at this stage, we define $\varphi_a$ by
\beq
\label{2.25e}
\partial_t \varphi_a = {\check B}_L (v_a)
\eeq

\noi with $\varphi_a (1) = 0$ and we define
\beq
\label{2.26e}
s_a = \nabla \varphi_a \quad , \quad B_a = B(v_a) \quad , \quad K_a = s_a + B_a \quad , \quad {\check B}_a = {\check B}(v_a)\ .
\eeq

\noi We next define the difference variables
\beq
\label{2.27e}
(w, \psi ) = \left ( v-v_a, \varphi - \varphi_a\right ) \ ,
\eeq
\beq
\label{2.28e}
G = B(v) - B_a \left ( =B(w, 2v_a + w)\right ) \ ,
\eeq
\beq
\label{2.29e}
\sigma = \nabla \psi \quad , \qquad L = \sigma + G
\eeq

\noi so that $K = K_a + L$. Substituting the definitions (\ref{2.27e})-(\ref{2.29e}) into the system (\ref{2.21e}) (\ref{2.22e}) yields the new system for $(w, \psi )$
\beq
\label{2.30e}
\partial_t \psi = {\check G}_L \left ( = {\check B}_L (w, 2v_a + w)\right )
\eeq
\beq
\label{2.31e}
i \partial_t w = Hw + H_1 v_a - R
\eeq

\noi where
\beq
\label{2.32e}
H_1 = i L \cdot \nabla_{K_a} + (i/2) \nabla \cdot \sigma + (1/2) L^2 + {\check G}_S + g(w, 2v_a + w) \ ,
\eeq
\beq
\label{2.33e}
R = i \partial_t v_a + (1/2) \Delta_{K_a} v_a - \left ( {\check B}_{aS} + g(v_a) \right ) v_a \ .
\eeq
\noi Again for technical reasons, it will be useful to consider also the partly linearized system for a new variable $w'$
\beq
\label{2.34e}
i \partial_t w' = H w' + H_1 v_a - R
\eeq

\noi where $H$ and $H_1$ are still associated with $(v, \varphi )$ (or $(w, \psi )$).\par

The remainder $R$ expresses the failure of $(v_a, \varphi_a)$ to satisfy the system (\ref{2.21e}) (\ref{2.22e}) and will have to tend to zero at a suitable rate in order to make it possible to solve that system.\par

The construction of solutions $(v, \varphi )$ of the system (\ref{2.21e}) (\ref{2.22e}) with prescribed asymptotic behaviour at $t= 0$ will be performed in two steps. The first step consists in solving the system (\ref{2.30e}) (\ref{2.31e}) with $(w, \psi )$ tending to zero as $t \to 0$ under assumptions on $(v_a, \varphi_a)$ of a general nature, the most important of which being decay assumptions on $R$ as $t \to 0$. This is done by first solving the linearized system (\ref{2.34e}) for $w'$, for given $(w, \psi)$ tending to zero as $t\to 0$, with $w'$ tending to zero as $t \to 0$. For that purpose, one first solves the Cauchy problem for the system  (\ref{2.34e}) with initial condition $w'(t_0) = 0$ for some $t_0 > 0$ and one takes the limit of the solution thereby obtained as $t_0 \to 0$. This procedure defines a map $\Gamma : w\to w'$. One then proves by a contraction method that the map $\Gamma$ has a fixed point in a suitable function space.\par

The second step of the method consists in choosing the asymptotic function $v_a$ so as to ensure the assumptions needed for the first step, and in particular the time decay of $R$. In the present problem, this will be simply achieved by taking $v_a = U(t) v_+$ for a suitably regular $v_+$. Substituting the previous results into the formulas (\ref{2.6e}) (\ref{2.16e}) will yield the corresponding results for the original system (\ref{1.1e}). In particular the solution $u$ thereby obtained will behave asymptotically as $u_a$ defined by 
\beq
\label{2.35e}
\widetilde{u}_a(t) = \overline{F \widetilde{u}_{ca} (1/t)}
\eeq
\beq
\label{2.36e}
u_{ca} = v_a \exp \left  (- i \varphi_a \right )
\eeq

\noi in analogy with (\ref{2.6e}) (\ref{2.16e}).\par

We now give a heuristic preview of the main results of this paper, stripped from most technicalities. They will be stated in full mathematical detail in Propositions 4.2-4.5 as regards the Cauchy problem with finite initial time and in Propositions 5.6, 5.7 as regards the Cauchy problem with initial time $t=0$ for $(v, \varphi )$ and $t = \infty$ for $u$. In order to state the results, we shall use the spaces $V^k$, $\Sigma^k$, $H_>^k$ and $H_>^{\infty}$ defined by (\ref{3.4e}) (\ref{3.9e}) (\ref{3.1e}) (\ref{3.2e}) below. In all those results we assume that $1 < k < 2$. The lower bound $k > 1$ plays an essential r\^ole, while the upper bound $k < 2$ is imposed only for convenience. It could be dispensed with at the expense of minor modifications of the proofs.The results for finite initial time can be summarized as follows.\\

\noi {\bf Proposition 2.1.} {\it Let $1 < k < 2$. \par

(1) Let $v_0 \in V^k$. For $\tau_0 > 0$, $\tau_0$ sufficiently small, there exists a unique solution $(v, \varphi )$ of the system (\ref{2.21e}) (\ref{2.22e}) such that $(v, \varphi ) (t_0) = (v_0, 0)$, $v \in ({\cal C} \cap L^{\infty})(I, V^k)$, $\varphi \in {\cal C} (I, H_>^{\infty})$ where $I = (0, \tau_0]$ and $(v, \varphi )$ is estimated in those spaces. Furthermore there exists $v_+ \in V^k$ such that $v(t$) tends to $v_+$ when $t \to 0$, and $(v, \varphi )$ behaves asymptotically as $(v_a, \varphi_a)$ when $t \to 0$, with $v_a = U(t) v_+$ and $\varphi_a$ a solution of (\ref{2.25e}). \par

(2) Let $\widetilde{u}_0 \in FV^k$. For $t_0$ sufficiently large, there exists a unique solution $u$ of the system (\ref{1.1e}) such that $u(t_0) = U(t_0) \widetilde{u}_0$, $\widetilde{u} \in {\cal C} (I, FV^k)$ where $I = [t_0, \infty )$, and $u$ is estimated in that space. Furthermore $u$ behaves as $u_a$ when $t \to \infty$, with $u_a$ defined by (\ref{2.35e}) (\ref{2.36e}) and $(v_a, \varphi_a)$ as in Part (1). \par

(3) Parts (1) and (2) hold with $V^k$ replaced everywhere by $\Sigma^k$. Furthermore for $u_0 \in \Sigma^k$, $u_0$ sufficiently small, there exists a unique solution $u \in {\cal C}({I\hskip-1truemm R} , \Sigma^k)$ of the system (\ref{1.1e}) with $u(0) = u_0$.} \\

We next summarize the results for zero or infinite initial time.\\

\noi {\bf Proposition 2.2.} {\it Let $1 < k < 2$. Let $v_+ \in V^{k+1}$, let $v_a = U(t) v_+$, let $\varphi_a$ be defined by (\ref{2.25e}) with $\varphi_a (1) = 0$ and let $u_a$ be defined by (\ref{2.35e}) (\ref{2.36e}). \par

(1) There exists $\tau > 0$ and there exists a unique solution $(v, \varphi )$ of the system (\ref{2.21e}) (\ref{2.22e}) such that $v \in ({\cal C} \cap L^{\infty}) (I, V^k)$ and $\varphi \in {\cal C}(I, H_>^{k+2})$, where $I = (0, \tau ]$, and such that $(v, \varphi )$ behaves asymptotically as $(v_a, \varphi_a)$ when $t \to 0$, in the sense that the difference $(v - v_a, \varphi - \varphi_a)$ tends to zero in suitable norms and at suitable rates when $t \to 0$. \par

(2) There exists $T > 0$ and there exists a unique solution $u$ of the system (\ref{1.1e}) such that $\widetilde{u} \in {\cal C} (I, FV^k)$, where $I = [T, \infty )$, and such that $u$ behaves asymptotically as $u_a$ when $t \to \infty$, in the sense that the difference $u - u_a$ tends to zero in suitable norms and at suitable rates when $t \to \infty$.}\\

\noi {\bf Remark 2.1.} There is a loss of one derivative from the asymptotic data $v_+$ to the solution $v$ in Proposition 2.2, part (1), so that Propositions 2.1 and 2.2 cannot be considered as the converse of each other. Furthermore the convergence properties of $(v, \varphi )$ to its asymptotic form $(v_a, \varphi_a)$ required in Proposition 2.2, part (1) to solve the system (\ref{2.21e}) (\ref{2.22e}) with initial time zero are stronger than those obtained for the solutions constructed in Proposition 2.1, part (1). A similar remark applies to the pair $(u, u_a)$.\\

We now describe the contents of the technical parts of this paper, namely Sections 3-5. In Section 3, we introduce some notation, we define the relevant function spaces and we collect a number of preliminary estimates. In Section 4 we study the Cauchy problem for finite initial time. We solve that problem for the auxiliary system (\ref{2.21e}) (\ref{2.22e}) (Proposition 4.2) and for the original system (\ref{1.1e}) (Proposition 4.4). We prove in particular the existence of small global solutions for the latter. We then analyse the asymptotic behaviour of the solutions thereby obtained for the auxiliary system (Proposition 4.3) and for the original system (Proposition 4.5). In Section 5 we study the Cauchy problem with initial time zero for $(v, \varphi )$ or infinity for $u$. We first give uniqueness results for $(v, \varphi )$ (Proposition 5.1) and for $u$ (Proposition 5.2). We then solve the Cauchy problem with prescribed asymptotic behaviour $(v_a, \varphi_a)$ as $t \to 0$ for $(v, \varphi )$ (Propositions 5.3-6) and with prescribed asymptotic behaviour $u_a$ as $t \to \infty$ for $u$ (Proposition 5.7).

\mysection{Notation and preliminary estimates}
\hspace*{\parindent} In this section we introduce some notation and we collect a number of estimates which will be used throughout this paper. We denote by $\parallel \ \parallel_r$ the norm in $L^r \equiv L^r ({I\hskip-1truemm R}^n)$, to be used mostly in ${I\hskip-1truemm R}^2$, and by $<,>$ the scalar product in $L^2$. We shall use the Sobolev spaces $H^k \equiv H^k ({I\hskip-1truemm R}^n)$ defined for $k \in {I\hskip-1truemm R}$ by
$$H^k = \left \{ u \in {\cal S} '({I\hskip-1truemm R}^n ): \parallel u; H^k \parallel\ = \ \parallel <\omega >^k u \parallel_2 \ < \infty \right \} \ ,$$

\noi where $<\cdot > = (1 + |\cdot |^2)^{1/2}$, and $\omega = (- \Delta )^{1/2}$. Besides the standard Sobolev spaces $H^k$, we will use the associated homogeneous spaces $\dot{H}^k$ with norm $\parallel u; \dot{H}^k\parallel\ =\break \noindent \parallel \omega^k u \parallel_2$. If $0 < k < n/2$ it is understood that $\dot{H}^k \subset L^r$ with $k = n/2 - n/r$. In addition we shall use the notation
\beq
\label{3.1e}
H_>^k = \ \mathrel{\mathop \cap_{0 < \ell \leq k}}\dot{H}^{\ell}
\eeq

\noi and
\beq
\label{3.2e}
H_>^{\infty} = \ \mathrel{\mathop \cap_{\ell > 0}}\dot{H}^{\ell}\  .
\eeq

\noi For any Banach space $X \subset {\cal S}' ({I\hskip-1truemm R}^n)$, we use the notation 
$$FX = \left \{ v \in  {\cal S}' ({I\hskip-1truemm R}^n) : F^* v \in X \right \}\ .$$

\noi For any $k \geq 0$, $\ell \geq 0$, we define the space
\beq
\label{3.3e}
H^{k, \ell} = \left \{ v \in {\cal S}' ({I\hskip-1truemm R}^n) : \ \parallel v; H^{k, \ell}\parallel\ = \ \parallel <x>^{\ell} <\omega >^k u \parallel_2 \ < \infty \right \}\ .
\eeq

\noi In particular $H^k = H^{k,0}$ and $F H^k = H^{0, k}$. For $1 < k  < 2$ we shall make extensive use of the space $V^k$ defined by
\beq
\label{3.4e}
V^k = \left \{ v \in {\cal S}' ({I\hskip-1truemm R}^n) : \ \parallel v; V^{k}\parallel\ = \ \parallel <\omega >^k v \parallel_2 \ \vee \ \parallel <\omega >^{k-1} x v \parallel_2 \ < \infty \right \}\ ,
\eeq

\noi where for real numbers $a$ and $b$ we use the notation $a \vee b = {\rm Max} (a,b)$ and $a\wedge b = {\rm Min} (a,b)$. Clearly $V^k = H^k \cap H^{k-1,1}$. More generally, for $0 \leq \rho \leq 1$ we define the spaces
\beq
\label{3.5e}
V^{k , \rho } = H^{k-1 + \rho } \cap H^{k-1, \rho}
\eeq

\noi so that
$$V^{k,0} = H^{k-1} \qquad , \qquad V^{k,1} = V^k \ .$$

\noi The spaces $V^{k, \rho }$ interpolate between $H^{k-1}$ and $V^k$. \par

From the commutation relation
\beq
\label{3.6e}
xU(-t) = U(-t) (x + it \nabla )
\eeq

\noi it follows that the space $V^k$ is invariant under the operator $U(t)$ and that 
\beq
\label{3.7e}
| \ \parallel v; V^k \parallel \ - \ \parallel U(t) v; V^k \parallel \ | \ \leq \ |t| \parallel v; H^k \parallel
\eeq

\noi so that
\beq
\label{3.8e}
\parallel U(t) v; V^k \parallel \ \leq (1 + |t|) \parallel v; V^k \parallel \ .
\eeq

\noi However the space $FV^k$ is not invariant under $U(t)$. For that and other reasons we shall also use the smaller space
\beq
\label{3.9e}
\Sigma^k = H^{k,0} \cap H^{0, k}
\eeq

\noi where $F$ acts an an isometry. From general interpolation theory it follows that $\Sigma^k \subset V^k$ with 
\beq
\label{3.10e}
\parallel v; V^k \parallel \ \leq \ C \parallel v; \Sigma^k\parallel
\eeq

\noi \cite{20newr} \cite{21r}. Furthermore $\Sigma^k$ is invariant under the evolution operator $U(t)$ and
\bea
\label{3.11e}
&&\parallel  U(-t) v; \Sigma^k \parallel \ = \ \parallel <\omega >^kv \parallel _2\ \vee \ \parallel <x+it\nabla >^k v \parallel _2\nn \\
&&\leq \ C \left ( \parallel  v ; \Sigma^k \parallel \ + \ |t|^k \parallel \omega^k v \parallel _2 \right ) \ \leq \ C (1 + |t|)^k \parallel  v; \Sigma^k\parallel \ .
\eea

\noi For the reader's convenience, we give a simple direct proof of those two facts and in particular of (\ref{3.10e}) (\ref{3.11e}) in the Appendix.\par

For any interval $I$ and for any Banach space $X$, we denote by ${\cal C}(I, X)$ (resp. ${\cal C}_w(I, X))$ the space of strongly (resp. weakly) continuous functions from $I$ to $X$ and by $L^{\infty} (I, X)$ (resp. $L_{loc}^{\infty} (I, X))$ the space of measurable essentially bounded (resp. locally essentially bounded) functions from $I$ to $X$. For $I$ an open interval, we denote by ${\cal D}' (I, X)$ the space of vector valued distributions form $I$ to $X$. We shall say that an evolution equation has a solution in $I$ with values in $X$ if the equation is satisfied in ${\cal D}'(I_0, X)$, where $I_0$ is the interior of $I$.\par

We shall use extensively the following Sobolev inequalities, stated here in ${I\hskip-1truemm R}^n$, but used only in ${I\hskip-1truemm R}^2$, and the following Leibnitz and commutator estimates.\\

\noi {\bf Lemma 3.1.} {\it (1) Let $1 < r \leq \infty$, $1 < r_1, r_2 < \infty$ and $0 \leq j < \ell $. If $r = \infty$, assume in addition that \ $\ell -j > n/r_2$. Let $\sigma$ satisfy \ $j/\ell  \leq \sigma \leq 1$ and
$$n/r - j = (1 - \sigma )n/r_1 + \sigma (n/r_2 - \ell )\ .$$

\noi Then the following estimate holds~:
$$\parallel \omega^j v \parallel_r\ \leq \ C \parallel v \parallel_{r_1}^{1-\sigma} \ \parallel \omega^{\ell} v \parallel_{r_2}^{\sigma}\ .$$

(2) Let $1 < r, r_1, r_3 < \infty$ and
$$1/r = 1/r_1 + 1/r_2 = 1/r_3 + 1/r_4\ .$$

\noi Then the following estimates hold~:
$$\parallel \omega^{\ell} (v_1 v_2) \parallel_r\ \leq \ C\left ( \parallel \omega^{\ell} v_1 \parallel_{r_1}\  \parallel v_2  \parallel_{r_2}\ + \  \parallel \omega^{\ell} v_2  \parallel_{r_3}\  \parallel v_1  \parallel_{r_4}\right )$$

\noi for $\ell \geq 0$, and 
$$ \parallel [\omega^{\ell} , v_1]v_2  \parallel_r\ \leq \ C\left ( \parallel \omega^{\ell} v_1 \parallel_{r_1}\  \parallel v_2  \parallel_{r_2}\ + \  \parallel \omega^{\ell -1} v_2  \parallel_{r_3}\  \parallel \nabla v_1  \parallel_{r_4}\right ) $$

\noi for $\ell  \geq 1$, where $[\ ,\ ]$ denotes the commutator.\par

In particular for $n=2$ and $0 < \ell <1$
\beq
\label{3.12e}
\parallel \omega^{\ell} (v_1 v_2) \parallel_2\ \leq \ C\parallel \omega^{\ell} v_1 \parallel_2  \left ( \parallel v_2  \parallel_{\infty}\ + \  \parallel \nabla v_2  \parallel_{2}\right ) \ .
\eeq

(3) Let $0 < \ell < 1$. Then the following estimate holds~:}
\beq
\label{3.13e}
\parallel [\omega^{\ell} ,v_1 ]\nabla v_2 \parallel_2\ \leq \ C\parallel F\nabla v_1 \parallel_1  \ \parallel \omega^{\ell} v_2   \parallel_2 \ .
\eeq 

\noi {\bf Proof.}  Part (1) follows from the Hardy-Littlewood-Sobolev (HLS) inequality \cite{19r} (from the Young inequality if $r = \infty$), from Paley-Littlewood theory and interpolation. Part (2) is proved in \cite{11r} \cite{13r} with $\omega$ replaced by $<\omega >$ and follows therefrom by a scaling argument.\par

Part (3) follows from the relation
$$\left ( F [\omega^{\ell} ,v_1 ]\nabla v_2\right ) (\xi ) = i \int d \eta \widehat{v}_1 (\xi - \eta ) \left ( |\xi |^{\ell} - |\eta |^{\ell}\right ) \eta \widehat{v}_2 (\eta ) \  ,$$

\noi from the inequality
$$|\ | \xi |^{\ell} - |\eta |^{\ell} | \ |\eta | \leq |\xi - \eta |\ |\eta |^{\ell}$$

\noi and from the Young inequality. \par \nobreak \hfill $\sq$ \par

We shall also use following lemma.\\

\noi {\bf Lemma 3.2.} {\it Let $1 < k < 2$. Then the following estimates hold~:}
\beq
\label{3.14e}
\parallel \omega^{2k-2} (xv_1 v_2) \parallel_2\ \leq \ C\parallel v_1 ; V^k  \parallel\   \parallel v_2 ; V^k   \parallel \ , 
\eeq
\beq
\label{3.15e}
\parallel \omega^{k+1} (xv_1 v_2) \parallel_2\ \leq \ C\parallel v_1 ; V^{k+1}  \parallel\   \parallel v_2 ; V^{k +1}  \parallel \ .
\eeq
\vskip 5 truemm

\noi {\bf Proof.} For $\ell \geq 1$ we write
$$\omega^{\ell} (xv_1 v_2) = [ \omega^{\ell} , v_1] x v_2 + v_1   [ \omega^{\ell} , x]v_2 + x v_1  \omega^{\ell} v_2  \ .$$

\noi By Lemma 3.1, part (2) we estimate
\begin{eqnarray*}
\parallel \omega^{\ell} (xv_1 v_2) \parallel_2 &\leq& C\Big ( \parallel \omega^{\ell} v_1  \parallel_{r_1} \   \parallel xv_2  \parallel_{r_2} \ + \  \parallel \nabla v_1  \parallel_{r_2}\ \parallel \omega^{\ell - 1} (xv_2)  \parallel_{r_1}\\
&+& \parallel v_1 \parallel_{r_2}\ \parallel \omega^{\ell - 1} v_2  \parallel_{r_1}\ +\ \parallel xv_1  \parallel_{r_2}\ \parallel \omega^{\ell } v_2  \parallel_{r_1}\Big )
\end{eqnarray*}

\noi with $1/r_1 + 1/r_2 = 1/2$, $2 \leq r_1 < \infty$. For $\ell = k+1$ we take $r_1 = 2$, $r_2 = \infty$ which immediately implies (\ref{3.15e}). For $k \geq 3/2$ and $\ell = 2k-2$ we take $2/r_1 = k-1$, $2/r_2 = 2 - k$ and we apply Lemma 3.1, part (1) to obtain 
\begin{eqnarray*}
\parallel \omega^{2k-2} (xv_1 v_2) \parallel_2 &\leq& C\Big ( \parallel \omega^{k} v_1  \parallel_2 \   \parallel \omega^{k-1} (xv_2)  \parallel_{2} \ + \  \parallel \omega^{k-1} (xv_1)  \parallel_{2}\ \parallel \omega^{k} v_2  \parallel_{2}\\
&+& \parallel \omega^{k - 1} v_1  \parallel_{2}\  \parallel \omega^{k-1} v_2  \parallel_{2}\Big )
\end{eqnarray*}

\noi which implies (\ref{3.14e}) in that case. For $k < 3/2$ so that $0 < \ell = 2 k-2 < 1$ we use the fact \cite{19r} that 
\beq
\label{3.16e}
\parallel \omega^{\ell} f\parallel_2^2\  = C \int dy |y|^{-2-2\ell} \parallel  \tau_y f - f \parallel_2^2
\eeq

\noi where $(\tau_y f)(x) = f(x+y)$. Taking $f(x) = xv_1(x) v_2(x)$ we write
\begin{eqnarray*}
\tau_y (xv_1 v_2) - x v_1 v_2 &=& \left ( \tau_y (xv_1)\right ) \left ( \tau_y v_2 - v_2\right ) + \left ( \tau_y v_1 - v_1\right ) xv_2\\ 
&&+ y \left ( \tau_y v_1\right ) v_2
\end{eqnarray*}

\noi so that from (\ref{3.11e}) we obtain
\bea
\label{3.17e}
&&\parallel \omega^{\ell} (xv_1 v_2) \parallel_2^2 \ \leq \ C\left ( \int dy|y|^{-2-2\ell} \Big ( \parallel \tau_y v_1  -v_1\parallel_{r_1}^2  \   \parallel xv_2  \parallel_{r_2}^2 \  \right .\nn \\
&&\left . + \  \parallel \tau_y v_2 - v_2  \parallel_{r_1}^2 \ \parallel xv_1  \parallel_{r_2}^2 \Big )+  \int dy |y|^{-2\ell} \parallel (\tau_y v_1)v_2 \parallel_{2}^2\right  )
\eea

\noi where $2/r_1 = k-1$, $2/r_2 = 2 - k$. We estimate
\begin{eqnarray*}
&&\int dy |y|^{-2-2\ell} \parallel \tau_y v_1 - v_1 \parallel_{r_1}^2\ \leq \ C \int dy |y|^{-2-2\ell} \parallel \omega^{2-k} \left ( \tau_y v_1 - v_1 \right ) \parallel_2^2\\
&&=\  C \parallel \omega^k v_1 \parallel_2^2
\end{eqnarray*}

\noi by Lemma 3.1, part (1), by (\ref{3.16e}) and the fact that $\omega$ commutes with $\tau_y$. We estimate $\parallel xv_2 \parallel_{r_2}$ by Lemma 3.1, part (1), we estimate the second term in (\ref{3.17e}) in the same way and we estimate the last term in (\ref{3.17e}) by the HLS inequality and Lemma 3.1, part (1) again. This yields (\ref{3.15e}). \par \nobreak \hfill $\sq$ \par

We next derive some estimates of the functions $B$ and $g$ defined by (\ref{2.11e}) (\ref{2.12e}) (\ref{2.13e}).\\

\noi {\bf Lemma 3.3.} {\it Let $1 < k < 2$. Then the following estimates hold~:} 
\beq
\label{3.18e}
\parallel \omega^{\ell} B(v_1 , v_2) \parallel_2\  \leq \ C \parallel v_1; H^k \parallel\ \parallel v_2 ; H^k \parallel\qquad {\it for}\ 0 < \ell \leq k+1\ ,
\eeq
\beq
\label{3.19e}
\parallel \omega^{\ell} B(v_1 , v_2) \parallel_2\  \leq \ C \parallel v_1 \parallel_2\ \parallel v_2 ; H^k \parallel\qquad {\it for}\ 0 < \ell \leq 1\ ,
\eeq
\beq
\label{3.20e}
\parallel \omega^{\ell} {\check B}(v_1 , v_2) \parallel_2\  \leq \ C\ t^{-1} \parallel v_1; V^k \parallel\ \parallel v_2 ; V^k \parallel\qquad {\it for}\ 0 < \ell \leq 2k-1\ ,
\eeq
\beq
\label{3.21e}
\parallel \omega^{\ell} {\check B}(v_1 , v_2) \parallel_2\  \leq \ C\ t^{-1} \parallel v_1 \parallel_2\ \parallel v_2 ; V^k \parallel\qquad {\it for}\ 0 < \ell \leq k-1\ ,
\eeq
\beq
\label{3.22e}
\parallel \omega^{\ell} {\check B}(v_1,  v_2) \parallel_2\  \leq \ C\ t^{-1} \parallel v_1;H^{k-1} \parallel\ \parallel v_2 ; V^k \parallel\qquad {\it for}\ 0 < \ell \leq 2(k-1)\ ,
\eeq
\beq
\label{3.23e}
\parallel \omega^{\ell}B(v)\parallel_2 \ \vee \ t\parallel \omega^{\ell}  {\check B}(v) \parallel_2\  \leq \ C\parallel v;V^{k+1} \parallel^2 \qquad {\it for}\ 0 < \ell \leq k+2\ ,
\eeq
\beq
\label{3.24e}
\parallel \omega^{\ell} g(v_1, v_2) \parallel_2\  \leq \ C \parallel v_1; H^k \parallel\ \parallel v_2 ; H^k \parallel\qquad {\it for}\ 0 \leq \ell \leq k\ ,
\eeq
\beq
\label{3.25e}
\parallel g(v_1, v_2) \parallel_2\  \leq \ C \parallel v_1 \parallel_2\ \parallel v_2 ; H^k \parallel \ . 
\eeq

\vskip 5 truemm

\noi {\bf Proof.} The estimates (\ref{3.18e}) (\ref{3.19e}) (\ref{3.24e}) (\ref{3.25e}) and the estimate of $B$ in (\ref{3.23e}) follow from Lemma 3.1 possibly supplemented by the HLS inequality for $B$ if $\ell < 1$. In order to derive the estimates for ${\check B}$ we first remark that 
\beq
\label{3.26e}
{\check B}(v_1, v_2) = 2t^{-1} \Delta^{-1} \partial_k \ {\rm Im} \left ( \overline{v}_{1k} \ x \cdot v_2 + \overline{v}_{2k} \ x \cdot v_1 \right ) \ .
\eeq

\noi This follows formally by commuting $x$ with $\Delta^{-1} \partial_k$ and can be proved by a regularization and a limiting procedure. The estimates (\ref{3.21e}) (\ref{3.22e}) use only the assumptions that $xv_2 \in H^{k-1}$ and follow from (\ref{3.26e}), from Lemma 3.1 and from the HLS inequality if $\ell < 1$. Finally  (\ref{3.20e}) and the estimate of ${\check B}$ in (\ref{3.23e}) follow from Lemma 3.2. \par \nobreak \hfill $\sq$ \par

We shall need also the following estimates.\\

\noi {\bf Lemma 3.4.} {\it Let $1 < k < 2$ and $k \leq m \leq 2$. Let $v \in V^k$ and $\nabla \varphi \in H^{m-1}$ with
\beq
\label{3.27e}
a = \ \parallel v; V^k\parallel \qquad , \quad \mu = \ \parallel \nabla \varphi ; H^{m-1} \parallel \ .
\eeq

(1) The following estimates hold~: 
\beq
\label{3.28e}
\parallel v \exp (- i \varphi ) \parallel_2\ = \ \parallel v \parallel_2 \ \leq a \ ,
\eeq
\beq
\label{3.29e}
\parallel xv \exp (- i \varphi ) \parallel_2\ = \ \parallel xv \parallel_2 \ \leq a \ ,
\eeq
\beq
\label{3.30e}
\parallel \omega^{k-1} xv \exp (- i \varphi ) \parallel_2\ \leq  C \ a(1 + \mu )  \ ,
\eeq
\beq
\label{3.31e}
\parallel \omega^{k} v \exp (- i \varphi ) \parallel_2\ \leq  C \ a(1 + \mu )^{1+(k-1)/(m-1)}  \ .
\eeq

(2) Let in addition $\varphi \in L^{\infty}$ with $\parallel \varphi \parallel_{\infty}\  \leq \mu$. Then 
\beq
\label{3.32e}
\parallel v \left ( \exp (- i \varphi ) - 1 \right )  \parallel_2\ \vee  \ \parallel xv \left ( \exp (- i \varphi ) - 1 \right ) \parallel_2 \ \leq a \mu \ ,
\eeq
\beq
\label{3.33e}
\parallel  \omega^{k-1}xv  \left ( \exp (- i \varphi ) - 1 \right )\parallel_2 \ \leq C \ a \ \mu \ ,
\eeq
\beq
\label{3.34e}
\parallel  \omega^{k} v  \left ( \exp (- i \varphi ) - 1 \right )\parallel_2\  \leq C \ a \ \mu (1 + \mu )^{(k-1)/(m-1)} \ .
\eeq

(3) Let in addition $xv \in L^{\infty} \cap \dot{H}^1$ with $\parallel xv \parallel_{\infty} \ \vee \ \parallel  \nabla x v \parallel _2 \ \leq a$. Then} 
\beq
\label{3.35e}
\parallel  \omega^{\ell} xv  \exp (- i \varphi ) \parallel_2\  \leq  a(1 + \mu )^{\ell} \qquad \hbox{for $0 \leq \ell \leq 1$}\ .
\eeq
\vskip 5 truemm

\noi {\bf Proof.} \underline{Part (1).} (\ref{3.28e}) and (\ref{3.29e}) are obvious. We next estimate by Lemma 3.1 
$$\parallel  \omega^{k-1} xv  \exp (- i \varphi ) \parallel_2\  \leq \ C \parallel  \omega^{k-1} xv   \parallel_2\left ( 1 + \ \parallel  \nabla \varphi \parallel_2\right )$$ 

\noi which implies (\ref{3.30e}), and 
\beq
\label{3.36e}
\parallel \omega^{k} v \exp (- i \varphi ) \parallel_2\ \leq  \ C  \left ( \parallel \omega^{k} v \parallel_2\ + \parallel v \parallel_{\infty} \ \parallel  \omega^k \exp (- i \varphi ) \parallel_2 \right ) \ .
\eeq

\noi We then interpolate
\beq
\label{3.37e}
\parallel \omega^{k} \exp (- i \varphi ) \parallel_2\ \leq  \  \parallel \nabla \varphi  \parallel_2^{(m-k)/(m-1)} \   \parallel  \omega^m \exp (- i \varphi ) \parallel_2^{(k-1)/(m-1)} 
\eeq

\noi and we estimate by Lemma 3.1 again
\beq
\label{3.38e}
\parallel \omega^{m}  \exp (- i \varphi ) \parallel_2\ =  \ \parallel \omega^{m-1} \nabla \varphi \exp (- i \varphi ) \parallel_2\ \leq \ C  \parallel \omega^m \varphi \parallel_{2} \left (1 + \  \parallel  \nabla \varphi  \parallel_2 \right ) 
\eeq

\noi which together with  (\ref{3.36e})  (\ref{3.37e}) implies  (\ref{3.31e}).\\

\noi \underline{Part (2).}  (\ref{3.32e}) is obvious. We next estimate by Lemma 3.1
$$\parallel \omega^{k-1} xv\left (  \exp (- i \varphi )-1 \right )  \parallel_2\ \leq \ C  \parallel \omega^{k-1} xv \parallel_{2} \left (\parallel   \varphi  \parallel_{\infty}\ + \  \parallel   \nabla \varphi  \parallel_2\right ) $$

\noi which implies (\ref{3.33e}), and
$$\parallel \omega^{k} v\left (  \exp (- i \varphi )-1 \right )  \parallel_2\ \leq \ C  \left ( \parallel \omega^{k} v \parallel_{2}\ \parallel   \varphi  \parallel_{\infty}\ + \  \parallel  v  \parallel_{\infty} \  \parallel    \omega^{k}  \exp (- i \varphi )  \parallel_2\right )  $$

\noi which together with (\ref{3.37e}) (\ref{3.38e}) implies (\ref{3.34e}). \\

\noi \underline{Part (3).} (\ref{3.35e}) is proved by interpolation between (\ref{3.29e}) and
$$\parallel \nabla (xv \exp (- i \varphi ) \parallel_2 \ \leq \ \parallel \nabla xv \parallel_2\ + \ \parallel xv \parallel_{\infty} \ \parallel \nabla \varphi \parallel_2 \ .$$
 \nobreak \hfill $\sq$\par

In order to take into account the time decay of the norms of $w$ as $t \to 0$, we introduce a function $h \in {\cal C} ((0, 1], {I\hskip-1truemm R}^+)$ such that the function $\overline{h}(t) = t^{-(3-k)/2} h(t)$ be non decreasing in $(0, 1]$ and satisfy 
\beq
\label{3.39e}
\int_0^t dt' \ t{'}^{-1} \ \overline{h}(t') \leq C\ \overline{h}(t)
\eeq

\noi for some $C > 0$ and all $t \in (0, 1]$. We shall use functions of the type 
\beq
\label{3.40e}
\overline{h}(t) = t^{\lambda} (1 - \ell n\ t)^{\mu}
\eeq

\noi with $\lambda > 0$ which clearly satisfy (\ref{3.39e}). Strictly speaking that function is increasing in $(0, 1]$ for $\lambda \geq \mu$ only, but not for $\lambda < \mu$. In the latter case, one can remedy that fact either by restricting oneself from the start to the smaller interval $(0, \tau ]$ with $\tau = \exp (1 - \mu /\lambda )$, or by replacing the previous function by 
$$\overline{h}(t) = \ \mathrel{\mathop {\rm Sup}_{0 < t' \leq t}}t{'}^{\lambda} (1 - \ell n\ t')^{\mu}$$

\noi which also satisfies (\ref{3.39e}). In what follows we shall freely use (\ref{3.40e}) and not mention that point any more.\par

For any interval $I \subset (0, 1]$, we define the space
\beq
\label{3.41e}
X(I) = \left \{ v : v \in {\cal C} (I, V^k) \quad \hbox{\rm and} \quad \parallel  v; X(I) \parallel \ = \ \mathrel{\mathop {\rm Sup}_{t \in I}}\ h(t)^{-1}\ \parallel  v(t) ; V^k \parallel  < \infty \right \}\ .
\eeq

We finally give estimates of the short and long range parts of ${\check B} = {\check B}(v)$ for time dependent $v$, namely
\beq
\label{3.42e}
\parallel \omega^m {\check B} _S\parallel _2 \ \leq \ t^{(p-m)/2} \parallel \omega^p {\check B}_S\parallel_2 \ \leq\ t^{(p-m)/2} \parallel \omega^p {\check B} \parallel _2
\eeq

\noi for $m \leq p$, and similarly
\beq
\label{3.43e}
\parallel \omega^m {\check B} _L\parallel _2 \ \leq \ \left ( 2t^{-1/2}\right )^{m-p} \parallel \omega^p {\check B}_L\parallel_2 \ \leq\ \ \left ( 2t^{-1/2}\right )^{m-p} \parallel \omega^p {\check B} \parallel _2
\eeq

\noi for $m \geq p$.

\mysection{The Cauchy problem with finite initial time}
\hspace*{\parindent} In this section we study the Cauchy problem with finite initial time in a neighborhood of zero for the auxiliary system (\ref{2.21e}) (\ref{2.22e}) and in a neighborhood of infinity for the original system (\ref{1.1e}). The main results are the existence of solutions defined down to $t= 0$ for $(v, \varphi )$, obtained in Proposition 4.2, and up to $t = \infty $ for $u$, derived therefrom in Proposition 4.4. As a by product we also obtain the existence of small global solutions for $u$. Furthermore we derive some results on the asymptotic behaviour of $(v, \varphi )$ as $t \to 0$ and of $u$ as $t \to \infty$, stated in Propositions 4.3 and 4.5 respectively. We treat the various questions in two types of function spaces. The largest convenient spaces are $V^k$ for $v$ and $\widetilde{v}$, and correspondingly $FV^k$ for $\widetilde{u}$. As mentioned in Section 3, those spaces have the drawback that $FV^k$ is not stable under the free Schr\"odinger evolution $U(t)$. The largest smaller spaces where stability holds are the spaces $\Sigma^k$ and we also state the various results specialized to those smaller spaces. One of the reasons for doing so is that the restriction to $\Sigma^k$ is necessary when dealing with the problem of small global solutions for $u$ (see Proposition 4.4, part (3)).\par

We begin this section by deriving some preliminary estimates for solutions of the partly linearized system  (\ref{2.21e}) (\ref{2.24e}). We recall that $s = \nabla \varphi$ and we use the notation $a_+ = a \vee 0$. \\

\noi {\bf Lemma 4.1.} {\it Let $1 < k < 2$ and $I \subset (0, 1]$. Let $v \in {\cal C} (I, V^k)$ and let 
$$y \equiv y(t) = \  \parallel v(t) ; V^k \parallel \ .$$

(1) Let $s \in {\cal C}(I, H^{k+1})$ and let $v' \in {\cal C}(I, V^k)$ be solution of (\ref{2.24e}). Then $v'$ satisfies the following estimates for all $t \in I$~:
\beq
\label{4.1e}
\parallel v'(t) \parallel_2\ = C\ ,
\eeq
\beq
\label{4.2e}
\left | \partial_t \parallel xv' \parallel_2 \right | \ \leq \ \parallel \nabla v' \parallel_2\ + \ \parallel s \parallel_2 \ \parallel v' \parallel_{\infty}\ + \ C \ y^2 \parallel v'\parallel_2\ ,
\eeq
$$\left | \partial_t \parallel \omega^k v' \parallel_2 \right | \ \leq \ C \Big \{ \left ( \parallel \nabla s\parallel_{\infty}\ + \ \parallel \omega^2 s \parallel_2 \ +\ \parallel s \parallel_{\infty}^2 \ + \ y^2 t^{k-2} + y^4 \right ) \parallel \omega^k v'\parallel_2$$
\beq
\label{4.3e}
+ \ \left ( \parallel \omega^{k+1}  s\parallel_2\ + \left (  \parallel \omega^k s \parallel_2 \ + y^2 \right  ) \left (  \parallel s \parallel_{\infty} \ + y^2\right  ) + y^2  t^{-1+(k-1)/2} + y^4 \right  ) \parallel  v'\parallel_{\infty} \Big \} \ ,
\eeq
$$\left | \partial_t \parallel \omega^{k-1} x v' \parallel_2 \right | \ \leq \ C \Big \{  \parallel \omega^k v'\parallel_2\ + \ \parallel \omega^{k-1} s \parallel_2 \left ( \parallel v' \parallel_{\infty} \ + \ \parallel \nabla v'\parallel_2\right )$$
$$+ \ y^2 \parallel \omega^{k-1}  v' \parallel_2\ + \left (  \parallel F\nabla  s\parallel_1\ + \  \parallel \omega^2 s \parallel_2 \ + \left (  \parallel \nabla s \parallel_2 \ + y^2\right  ) \left (  \parallel  s \parallel_{\infty}\ + y^2 \right )\right . $$
\beq
\label{4.4e}
\left . + \ y^2 t^{k-2} + y^4\right ) \parallel \omega^{k-1}  xv' \parallel_2 \Big \} \  .
\eeq

\noi Let in addition $v' \in {\cal C}(I, \Sigma^k)$. Then $v'$ satisfies the following estimate for all $t \in I$
\bea
\label{4.5e}
&&\left | \partial_t \parallel <x>^k v' \parallel_2 \right | \ \leq \ k\left ( \parallel <x>^{k-1} \nabla v' \parallel_2\ \right . \nn \\
&&\left . + \left ( \parallel s \parallel_{\infty} \ + Cy^2 \right ) \parallel <x>^{k-1} v' \parallel_2\ + \ \parallel v' \parallel_2\right ) \ .
\eea

(2) Let $\varphi$ satisfy (\ref{2.21e}). Then $\varphi$ satisfies the estimate
\beq
\label{4.6e}
\parallel  \omega^{\ell} \partial_t \varphi \parallel_2\ \leq C\ y^2\ t^{-1-(\ell /2 + 1/2 - k)_+}
\eeq

\noi for all $\ell > 0$ and all $t \in I$.}\\

\noi {\bf Proof.} \underline{Part (1).} (\ref{4.1e}) is obvious. (\ref{4.2e}) follows immediately from (\ref{2.24e}), from the commutation relation
$$[x, H] = \nabla_K$$

\noi and from Lemma 3.3. We next estimate $\parallel \omega^k v' \parallel_2$. By standard energy methods followed by Lemma 3.1, we estimate
$$\left | \partial_t \parallel \omega^{k}  v' \parallel_2 \right | \ \leq \  \parallel [\omega^k , s+B]\nabla v'\parallel_2\ + \ \parallel \omega^{k} (\nabla \cdot s) v'\parallel_2 \ + \  \parallel \omega^k (s+B)^2v' \parallel_2$$
$$+ \ \parallel \omega^{k}  {\check B}_Sv' \parallel_2\ +\   \parallel  \omega^k gv' \parallel_2$$
$$\ \leq \ C \Big \{  \Big (  \parallel \nabla (s + B) \parallel_{\infty} \ + \ \parallel \omega^2 (s+B) \parallel_2\ + \ \parallel \nabla \cdot  s \parallel_{\infty} \ +  \  \parallel s + B \parallel_{\infty}^2$$
$$+ \ \parallel  {\check B}_S \parallel_{\infty} \ +\   \parallel  g\parallel_{\infty}\Big ) \parallel \omega^k v' \parallel_2 \ + \ \Big ( \parallel  \omega^{k+1} s\parallel _2 \ + \ \parallel \omega^k (s+B) \parallel_2 \ \parallel s+B\parallel_{\infty} $$
\beq
\label{4.7e}
+\ \parallel \omega^k {\check B}_S \parallel_2 +\   \parallel  \omega^k g \parallel_2 \Big ) \parallel v' \parallel_{\infty}  \Big \}
\eeq

\noi from which (\ref{4.3e}) follows by (\ref{3.42e}) and Lemma 3.3.\par

We next estimate $\parallel \omega^{k-1} xv'\parallel_2$. By standard energy methods followed by Lemma 3.1, we estimate 
$$\left | \partial_t \parallel \omega^{k-1}  xv' \parallel_2 \right | \ \leq \  \parallel \omega^{k-1} \nabla_{s+B} v'\parallel_2\ + \ \parallel [\omega^{k-1}, s+B] \nabla x v'\parallel_2 $$
$$+ \  \parallel \omega^{k-1} (\nabla \cdot s) xv' \parallel_2+ \ \parallel \omega^{k-1} (s+B)^2 xv' \parallel_2\ + \   \parallel \omega^{k-1} {\check B}_Sxv' \parallel_2$$
$$ +\   \parallel  \omega^{k-1}  gxv' \parallel_2$$
$$\ \leq \ C \Big \{   \parallel \omega^k v'\parallel_2 \ + \ \parallel \omega^{k-1} s \parallel_2 \left (  \parallel v'\parallel_{\infty} \ + \ \parallel \nabla v'\parallel_2\right  ) + \ \left ( \parallel B \parallel_{\infty} \ +  \  \parallel \nabla  B \parallel_2\right ) $$
$$\times \ \parallel \omega^{k-1} v' \parallel_2\ + \Big ( \parallel F \nabla  (s+B) \parallel_1\ + \  \parallel \nabla \cdot s \parallel_{\infty}\ + \ \parallel \nabla  \nabla  \cdot s \parallel_2 \ + \ \parallel s+B\parallel_{\infty}^2$$
$$+ \ \parallel s+B\parallel_{\infty}\ \parallel \nabla  (s+ B) \parallel_2\ + \  \parallel  {\check B}_S \parallel_{\infty} \ +\   \parallel  \nabla {\check B}_S \parallel_2\ +\ \parallel g \parallel_{\infty}\ + \ \parallel  \nabla g\parallel_2\Big )$$
\beq
\label{4.8e}
\times\ \parallel\omega^{k-1} xv' \parallel_2  \Big \}
\eeq

\noi from which (\ref{4.4e}) follows by (\ref{3.42e}) and Lemma 3.3. Finally from the commutation relation
\beq
\label{4.9e}
\left [ <x>^k, H\right ] = \left ( \nabla <x>^k \right ) \cdot \nabla_K + (1/2) \left ( \Delta <x>^k \right ) 
\eeq

\noi with
$$\nabla <x>^k = k <x>^{k-2}x \quad , \quad \Delta <x>^k = k^2 <x>^{k-2} - k (k-2) <x>^{k-4}\ ,$$

\noi we obtain
$$\left | \partial_t \parallel <x>^k v' \parallel_2 \right | \ \leq \  \parallel [<x>^k, H]v' \parallel_2$$
$$\leq k \left ( \parallel <x>^{k-1} \nabla v' \parallel_2\ + \ \parallel <x>^{k-1} Kv' \parallel_2\ + \ \parallel v' \parallel_2\right )$$

\noi from which (\ref{4.5e}) follows.\\

\noi \underline{Part (2)} follows immediately from (\ref{3.43e}) and from Lemma 3.3. \par \nobreak \hfill $\sq$ \par

We next derive some estimates of the difference of two solutions of (\ref{2.21e}) (\ref{2.24e}) associated with two different $v's$. We shall use the following notation. Let $f_i$, $i = 1,2$ be two functions or operators. We define $f_{\pm} = (1/2) (f_1 \pm f_2)$ so that $f_1 = f_+ + f_-$, $f_2 = f_+ - f_-$ and $(fg)_{\pm} = f_+ g_{\pm} + f_- g_{\mp}$. Let now $v'_i$, $i = 1,2$ be a pair of solutions of (\ref{2.24e}) associated with a pair $(v_i, s_i)$, $i = 1,2$. Then $v'_-$ satisfies the equation 
\beq
\label{4.10e}
i \partial_t v'_- = H_+ v'_- + H_- v'_+
\eeq

\noi where 
\beq
\label{4.11e}
H_+ = - (1/2) \Delta _{K_+} + (1/2) K_-^2 + {\check B}_{S_+} + g_+ \  ,
\eeq
\beq
\label{4.12e}
H_- = iK_- \cdot  \nabla _{K_+} + (i/2) \nabla \cdot K_- + {\check B}_{S_-} + g_- \  .
\eeq

We can now state the difference estimates of two solutions of (\ref{2.21e}) (\ref{2.24e}).\\

\noi {\bf Lemma 4.2.} {\it Let $1 < k < 2$ and $I \subset (0, 1]$. Let $v_i \in {\cal C}(I, V^k)$, $i = 1,2$ and let
\beq
\label{4.13e}
y = y(t) = \ \mathrel{\mathop {\rm Max}_{i=1,2}}\ \parallel v_i(t) ; V^k \parallel \ . 
\eeq

(1) Let $s_i \in {\cal C}(I, H^{k+1})$, $i = 1,2$, and let $v'_i \in {\cal C}(I, V^k)$, $i = 1,2$, be solutions of (\ref{2.24e}) associated with $(v_i, s_i)$. Then the following estimate holds for all $t\in I$~: 
$$\left | \partial_t \parallel v{'}_- \parallel_2 \right | \ \leq \ C \Big \{  \left  ( \parallel \omega^{2-k} s_-\parallel_2\ + \ y\parallel v_- \parallel_2\right )  \parallel \omega^k v'_+ \parallel_2$$
$$+ \Big (  \parallel s_- \parallel_2 \left (  \parallel s_+\parallel_{\infty}  + y^2 \right ) \ + \ y \parallel v_-\parallel_2 \left (  \parallel \omega^{k-1} s_+ \parallel_2\ + y^2 \right )$$
\beq
\label{4.14e}
+ \  \parallel \nabla \cdot s_- \parallel_2 \ + \ y \parallel  v_- \parallel_2\ t^{-1+(k-1)/2} \Big ) \parallel v'_+ \parallel_{\infty}\Big \} \ .
\eeq

(2) Let $\varphi_i$, $i = 1,2$, satisfy (\ref{2.21e}) with $v = v_i$. Then the following estimate holds for all $\ell > 0$ and for all $t \in I$~:}
\beq
\label{4.15e}
\parallel \omega^{\ell } \partial_t \varphi_-\parallel_2\ \leq \ C\ y \parallel v_- \parallel_2\ t^{-1-(1/2)(\ell + 1 - k)_+} \ .
\eeq
\vskip 5 truemm

\noi {\bf Proof.} \underline{Part (1).} From (\ref{4.10e}) we estimate
$$\left | \partial_t \parallel v{'}_- \parallel_2 \right | \ \leq \  \parallel H_-v'_+\parallel_2\ \leq \ \parallel K_- \parallel_{r_1}  \ \parallel \nabla v'_+\parallel_{r_2}$$
$$+ \Big (  \parallel s_- \parallel_2  \ \parallel K_+\parallel_{\infty}  \ + \ \parallel B_-\parallel_{r_1}\   \parallel K_+ \parallel_{r_2}$$
$$+ \ \parallel \nabla \cdot s_-\parallel_2 \ + \ \parallel {\check B}_{S-} \parallel_2 \ + \ \parallel g_-\parallel_2\Big ) \parallel v'_+ \parallel_{\infty}$$

\noi with $r_1 = 2/(k-1)$, $r_2 = 2/(2-k)$, from which (\ref{4.14e}) follows by (\ref{3.42e}) and Lemma 3.3 which implies in particular
\bea
\label{4.16e}
&&\parallel K_- \parallel_{r_1}  \ \leq\ C \parallel \omega^{2-k} K_-\parallel_2\ \leq \ C \left ( \parallel \omega^{2-k} s_-\parallel_2\ + \parallel v_-\parallel_2\ \parallel v_+;V^k\parallel\right ) \ ,\nn \\
&&\parallel \omega^{k-1}  {\check B}_-\parallel_2\ \leq \ C\ t^{-1}  \parallel v_-\parallel_2\ \parallel v_+;V^k\parallel \ .
\eea
\vskip 5 truemm

\noi \underline{Part (2).}  (\ref{4.15e}) follows immediately from  (\ref{3.43e}) and  (\ref{4.16e}).\par \nobreak \hfill $\sq$ \par

We now turn to the study of the Cauchy problem with finite initial time for the auxiliary system  (\ref{2.21e}) (\ref{2.22e}). The first step is to solve that problem for the linearized system (\ref{2.24e}).\\

\noi {\bf Proposition 4.1.} {\it Let $1 < k < 2$ and $I \subset (0, 1]$, let $t_0 \in I$ and $v'_0 \in V^k$. Let $v \in {\cal C} (I, V^k)$ and $s \in {\cal C}(I, H^{k+1})$. Then there exists a unique solution $v' \in {\cal C} (I, V^k)$ of the system (\ref{2.24e}) with $v'(t_0) = v'_0$. That solution satisfies the estimates (\ref{4.1e})-(\ref{4.4e}) of Lemma 4.1, part (1). The difference of two such solutions satisfies the estimate (\ref{4.14e}) of Lemma 4.2, part (1). Uniqueness actually holds in $L^{\infty} (I, V^k)$. If in addition $v'_0 \in \Sigma^k$, then $v' \in {\cal C}(I, \Sigma^k )$ and $v'$ satisfies the estimate (\ref{4.5e}). }\\

That proposition can be proved for instance by a parabolic regularization followed by a limiting procedure. A similar result in a more complicated context appears in Proposition 4.1 of \cite{4r}.\par

We now turn to the main technical result of this section, namely the existence of solutions of the auxiliary system (\ref{2.21e}) (\ref{2.22e}) with sufficiently small initial time, defined down to time zero.\\

\noi {\bf Proposition 4.2.} {\it Let $1 < k < 2$ and let $v_0 \in V^k$ with $\parallel  v_0;V^k \parallel \ = a$. Then \par

(1) There exists $\overline{\tau}_0$, $0 < \overline{\tau}_0 \leq 1$, such that for any $\tau_0$, $0 < \tau_0 \leq \overline{\tau}_0$, there exists a unique solution $(v, \varphi )$ of the system (\ref{2.21e}) (\ref{2.22e}) with $(v, \varphi ) (\tau_0) = (v_0, 0)$ and such that $v \in ({\cal C} \cap L^{\infty}) (I, V^k)$, where $I = (0, \tau_0]$. Furthermore $\varphi \in {\cal C} (I, H_>^{\infty})$ and $(v, \varphi )$ satisfy 
\beq
\label{4.17e}
\parallel  v; L^{\infty}(I, V^k)\parallel \ \leq C\ a \ ,
\eeq
\beq
\label{4.18e}
\parallel \omega^{\ell} \varphi\parallel _2\ \leq C\ a^2 \left ( t^{-(\ell /2 + 1/2 - k)_+} - \ell n\ t \right )
\eeq

\noi for all $\ell > 0$ and all $t \in I$. The time $\overline{\tau}_0$ depends on $a$ according to 
\beq
\label{4.19e}
a\ \overline{\tau}_0^{(k-1)/4}\ \leq C\ .
\eeq

In particular one can take $\overline{\tau}_0 = 1$ for $a$ sufficiently small. Uniqueness holds actually under the condition $v \in {\cal C}(I, V^k)$.\par

If in addition $v_0 \in \Sigma^k$, then $v \in ({\cal C} \cap L^{\infty})(I, \Sigma^k)$ and $v$ satisfies the estimate
\beq
\label{4.20e}
\parallel  v; L^{\infty}(I, \Sigma^k) \parallel \ \leq \ C \parallel  v_0; \Sigma^k \parallel \ .
\eeq

(2) The map $v_0 \to (v, \varphi )$ is continuous for fixed $\tau_0$ on the bounded sets of $V^k$, from the $L^2$-norm of $v_0$ to the norm of $(v, \varphi )$ in $L^{\infty}(J, H^{k'} \oplus H_>^{\ell})$ for any $k'$, $0 \leq k' < k$ and any $\ell > 0$, and in the weak $\star$ sense in $L^{\infty}(J, V^{k} \oplus H_>^{\ell})$ for any interval $J \subset\subset I$. \par

If in addition $v_0 \in \Sigma^k$, the continuity of $v$ extends to norm continuity in $L^{\infty}(J, \Sigma^{k'})$ for any $k'$, $0 \leq k' < k$, and to weak $\star$ continuity in  $L^{\infty}(J,  \Sigma^{k})$ for any interval $J \subset\subset I$.}\\

\noi {\bf Proof.} \underline{Part (1).} The proof consists in showing that the map $\Gamma : v \to v'$ defined by Proposition 4.1 with $v' (\tau_0) = v(\tau_0) = v_0$ and $s (\tau_0) = 0$ is a contraction on a suitable bounded set ${\cal R}$ of $({\cal C} \cap L^{\infty})(I, V^k)$ for the norm in $L^{\infty} (I, L^2)$, where $I = (0, \tau_0]$. We define
\beq
\label{4.21e}
{\cal R} = \left \{ v \in ({\cal C} \cap L^{\infty}) (I, V^k) : v(t_0) = v_0, \ \parallel v ; L^{\infty} (I, V^k) \parallel \ \leq Y \right \}
\eeq

\noi for some constant $Y$ to be chosen later.\par

We first show that ${\cal R}$ is mapped into itself by $\Gamma$. Integrating (\ref{4.6e}) between $\tau_0$ and $t$, we  obtain 
\beq
\label{4.22e}
 \parallel  \omega^{\ell} s  \parallel _2 \ \leq C\ Y^2 \left ( t^{-(\ell /2 + 1 - k)_+} - \ell n\ t \right )
\eeq

\noi for all $\ell \geq 0$ and all $t \in I$, so that in particular
$$ \parallel s \parallel_{\infty} \ \leq C\ Y^2 \left ( t^{-(3/2-k)_+} - \ell n\ t \right ) \ ,$$
$$ \parallel \nabla s \parallel_{\infty} \ \leq\   \parallel F \nabla s  \parallel_1 \ \leq C\ Y^2 t^{k-2}  \ .$$

\noi Let now
\beq
\label{4.23e}
y' \equiv y'(t) \equiv \  \parallel v'(t) ; V^k  \parallel \ .
\eeq

\noi Substituting the previous estimates into (\ref{4.2e})-(\ref{4.4e}), we obtain
\beq
\label{4.24e}
\left | \partial_t \parallel xv' \parallel_2 \right | \ \leq C \left ( 1 +  Y^2(1 - \ell n\ t)\right ) y' \ ,
\eeq 
\beq
\label{4.25e}
\left | \partial_t \parallel \omega^{k}  v' \parallel_2 \right | \ \leq \ C\left ( Y^2\ t^{-1 + (k-1)/2} + Y^4\ t^{-1+k/2}\left ( t^{-(3/2-k)_+} - \ell n \ t \right ) \right ) y'
\eeq
\beq
\label{4.26e}
\left | \partial_t \parallel \omega^{k-1}  xv' \parallel_2 \right | \ \leq \ C\left ( 1 + Y^2\ t^{k-2} + Y^4\left (  t^{-(3-2k)_+} + (\ell n\ t)^2  \right ) \right ) y'\ . 
\eeq

\noi Integrating (\ref{4.24e})-(\ref{4.26e}) over time yields
\beq
\label{4.27e}
y'(t) \leq a\exp \left \{ C \left ( \tau_0 + Y^2 \tau_0^{(k-1)/2} + Y^4 \tau_0^{k/2} \left ( \tau_0^{-(3/2-k)_+} - \ell n\ \tau_0 \right ) \right ) \right \}
\eeq

\noi so that by choosing $Y = Ca$, taking $\tau_0$ sufficiently small according to 
\beq
\label{4.28e}
a^2\ \tau_0^{(k-1)/2} \leq C
\eeq

\noi for suitable constants $C$ and using the fact that
\beq
\label{4.29e}
k/2 - (3/2 - k)_+ > k-1
\eeq

\noi we obtain 
\beq
\label{4.30e}
Y' \equiv \ \parallel y'; L^{\infty}(I) \parallel \ \leq Y \ .
\eeq

\noi This proves that ${\cal R}$ defined by (\ref{4.21e}) is mapped into itself by $\Gamma$.\par

We next prove that $\Gamma$ is a contraction for the norm in $L^{\infty} (I, L^2)$ on ${\cal R}$. Let $v_i \in {\cal R}$ and $v'_i = \Gamma v_i$, $i = 1,2$. We estimate the difference $v'_-$ by Lemma 4.2. From (\ref{4.15e}), we obtain 
\beq
\label{4.31e}
 \parallel  \omega^{\ell} s_-  \parallel _2 \ \leq C\ Y\ Y_- \ t^{- 1 +( k-\ell)/2}
\eeq

\noi for all $\ell \geq 0$ and all $t \in I$, where 
$$Y_- \equiv \  \parallel v_-;L^{\infty} (I, L^2)  \parallel \  .$$

\noi Substituting (\ref{4.31e}) into (\ref{4.14e}) and using the fact that $v'_+ \in {\cal R}$, we obtain 
$$\left | \partial_t \parallel v'_- \parallel_2 \right | \ \leq \ C\ Y_-\left ( Y^2\ t^{-1+ (k-1)/2} + Y^4\ t^{-1 + k/2} \left ( t^{-(3/2-k)_+} - \ell n\ t  \right ) \right )  $$

\noi so that by integration over time 
$$Y'_- \equiv\ \parallel v'_-; L^{\infty}(I, L^2) \parallel\ \leq C\ Y_- \left ( Y^2\ \tau_0^{(k-1)/2} + Y^4\ \tau_0^{k/2} \left ( \tau_0^{-(3/2-k)_+} - \ell n \ \tau_0 \right ) \right ) \ .$$

\noi Taking again $\tau_0$ sufficiently small according to (\ref{4.28e}) (possibly with a smaller constant) and using again (\ref{4.29e}), we obtain
$$Y'_- \leq (1/2)Y_-$$

\noi which proves that $\Gamma$ is a contraction for the $L^{\infty}(I, L^2)$ norm on ${\cal R}$. Together with the fact that ${\cal R}$ is closed for that norm, this proves that $\Gamma$ has a unique fixed point in ${\cal R}$. Uniqueness for $v \in {\cal C}(I, V^k)$ follows from similar estimates.\par

The last statement follows by integration of (\ref{4.5e}), using the fact that 
$$\parallel <x>^{k-1} \nabla v \parallel_2\ \leq\ C \parallel v; \Sigma^k\parallel \ .$$
\vskip 5 truemm

\noi \underline{Part (2).} Let $v_i$, $i = 1,2$, be two solutions of the previous type of the system (\ref{2.21e}) (\ref{2.22e}) with different initial data $v_{0i}$, $i = 1,2$. We estimate $v_-$ in $L^{\infty} (I, L^2)$ by using again Lemma 4.2, where now $v'_i = v_i$, but $v_- (\tau_0) = v_{0-} = (1/2)(v_{01} - v_{02})\not= 0$. By the same computation as in the contraction proof, we obtain 
$$Y_- \leq \ 2  \parallel v_{0-}\parallel_2$$

\noi which proves the continuity of $v$ from $L^2$ to $L^{\infty}(I,L^2)$. The remaining continuities follow therefrom by interpolation with boundedness of $v$ in $L^{\infty} (I, V^k)$ or in $L^{\infty}(I , \Sigma^k)$, from standard compactness arguments and from (\ref{4.15e}). \par \nobreak \hfill $\sq$ \par

\noi {\bf Remark 4.1.} We have considered the Cauchy problem for the system (\ref{2.21e}) (\ref{2.22e}) with initial condition $\varphi (\tau_0) = 0$. One could easily take instead an initial condition $\varphi (\tau_0) = \varphi_0$ for some $\varphi_0 \in H_>^{k+2}$ satisfying (\ref{4.18e}) for $0 < \ell \leq k+2$ and $t = \tau_0$. The solutions thereby obtained would exhibit continuity properties with respect to $\varphi_0$ and $\tau_0$.\\

We next derive asymptotic properties in time of the solutions of the auxiliary system (\ref{2.21e}) (\ref{2.22e}) obtained in Proposition 4.2. We prove in particular the existence of a limit $v_+$ of $v(t)$ when $t \to 0$ (the subscript $+$ is used here with a different meaning from that used when comparing two solutions) and we provide two asymptotic forms of the phase $\varphi$. The first one is more accurate while the second one has a simpler form.\\

\noi {\bf Propostion 4.3.} {\it Let $1 < k < 2$. Let $(v, \varphi )$ be a solution of the system (\ref{2.21e}) (\ref{2.22e}) as obtained in Proposition 4.2 and let
\beq
\label{4.32e}
Y = \ \parallel v;L^{\infty} (I, V^k)  \parallel
\eeq

\noi where $I = (0, \tau_0 ]$. Then \par

(1) There exists $v_+ \in V^k$ such that $v(t)$ tends to $v_+$ when $t \to 0$, strongly in $V^{k, \rho }$ for $0 \leq \rho < 1$ (and in particular in $H^{k-1}$) and weakly in $V^k$. Furthermore
\beq
\label{4.33e}
\parallel v_+ ; V^k \parallel\ \leq \ \mathrel{\mathop {\rm lim\ inf }_{t \to 0}}\ \parallel v(t) ; V^k \parallel\ \leq Y
\eeq

\noi and the following estimate holds for all $t\in I$
\beq
\label{4.34e}
\parallel  v(t) - v_a (t) ; H^{k-1}\parallel\ \leq C\ Y^3 (1 + Y^2) t^{k/2}
\eeq

\noi where
\beq
\label{4.35e}
v_a(t) = U(t) v_+ \ .
\eeq

Similar estimates hold in $V^{k, \rho}$, $0 \leq \rho < 1$, by interpolation between (\ref{4.32e}) (\ref{4.33e}) (\ref{4.34e}).\par

If in addition $v \in ({\cal C} \cap L^{\infty}) (I, \Sigma^k)$, then $v_+ \in \Sigma^k$ and $v_+$ satisfies 
\beq
\label{4.36e}
\parallel  v_+ ; \Sigma^k \parallel \ \leq \ \mathrel{\mathop {\rm lim\ inf }_{t \to 0}}\ \parallel v(t) ; \Sigma^k \parallel\ .
\eeq

\noi Furthermore $v(t)$ tends to $v_+$ when $t \to 0$ strongly in $\Sigma^{k'}$ for $0 \leq k' < k$ and weakly in $\Sigma^k$.\par

(2) Define $\varphi_1(t)$ by 
\beq
\label{4.37e}
\partial_t \varphi_1 = {\check B}_L (v_a) \quad , \qquad \varphi_1(1) = 0 \ .
\eeq

\noi Then there exists $\psi_{1+}\in \dot{H}^{\ell}$ for $0 < \ell < 3k-2$ such that $\varphi (t) - \varphi_1 (t)$ tends to $\psi_{1+}$ in $\dot{H}^{\ell}$ for all such $\ell$ when $t \to 0$. Define $\varphi_a = \varphi_1 + \psi_{1+}$. Then the following estimates hold for all such $\ell$ and all $t \in I$~: 
\beq
\label{4.38e}
\parallel   \omega^{\ell} \left ( \varphi (t) - \varphi_a(t) \right ) \parallel_2\ \leq C\ Y^4(1 + Y^2) t^{k/2 - (\ell /2 + 1 - k)_+}\ ,
\eeq
\beq
\label{4.39e}
\parallel   \omega^{\ell} \varphi_a(t) \parallel_2\ \leq C\ Y^2(1 + Y^2)\left (  t^{-(\ell /2 + 1/2 -k)_+} - \ell n\ t \right )  \  .
\eeq

(3) Define 
\beq
\label{4.40e}
\varphi_2(t) = (\ell n \ t) x \cdot B(v_+) \ .
\eeq

\noi Then there exists $\psi_{2+} \in \dot{H}^{\ell}$ for $0 < \ell < 2k-1$ such that $\varphi (t) - \varphi_2(t)$ tends to $\psi_{2+}$ in $\dot{H}^{\ell}$ for all such $\ell$ when $t \to 0$. Define $\varphi_b = \varphi_2 + \psi_{2+}$. Then the following estimates hold for all such $\ell$ and all $t \in I$~:}
\beq
\label{4.41e}
\parallel   \omega^{\ell} \left ( \varphi (t) - \varphi_b(t) \right ) \parallel_2\ \leq C\ Y^2(1 + Y^2)^2\  t^{k/2 -(1/2)  (\ell  + 1 - k)_+}\ ,
\eeq
\beq
\label{4.42e}
\parallel   \omega^{\ell} \varphi_b(t) \parallel_2\ \leq C\ Y^2(1 + Y^2)^2 \left (  t^{-(\ell /2 + 1/2 -k)_+} - \ell n\ t \right )  \  .
\eeq
\vskip 5 truemm

\noi {\bf Proof.} \underline{Part (1).} Let $\widetilde{v} = U(-t)v$. Then $\widetilde{v}$ satisfies the equation 
\beq
\label{4.43e}
i \partial_t \widetilde{v} = U(-t) \left \{ i (s+ B) \cdot \nabla v + \left ( (i/2) \nabla \cdot s + (1/2) (s+ B)^2 +  {\check B}_S + g \right ) v \right \} \ .
\eeq

\noi By Lemma 3.1, Lemma 3.3, (\ref{3.42e}) and (\ref{4.18e}), we estimate 
$$\parallel   \omega^{\ell} \partial_t \widetilde{v}\parallel _2\ \leq\ C \Big \{ \left  ( \parallel  \omega^{\ell +2 - k}(s+B)\parallel _2 \ + \  \parallel  \omega^{\ell +1 - k}(s+B)\parallel _{\infty}\right )  \parallel  \omega^k v \parallel_2$$ 
$$+ \left ( \parallel \omega^{\ell} \nabla \cdot s \parallel_2 \ + \ \parallel \omega^{\ell}(s + B)^2 \parallel_2\ + \ \parallel \omega^{\ell} {\check B}_S\parallel_2 \ + \ \parallel \omega^{\ell} g \parallel_2\right ) \left ( \parallel v \parallel_{\infty}\ + \ \parallel \nabla v  \parallel_2\right ) \Big \}$$
$$\leq C \left \{ Y^3\left ( t^{-(\ell /2 + 3/2 - k)_+} - \ell n \ t \right ) + Y^5 \left ( t^{-(\ell /2 + 5/2 - 2k)_+} + (\ell n\ t)^2 \right ) \right \}$$
\beq
\label{4.44e}
\leq C\ Y^3 (1 + Y^2) \left ( t^{-(\ell /2 + 3/2 - k)_+} + (\ell n \ t)^2 \right ) \leq C\  Y^3 (1 + Y^2)  t^{-1 + k/2} 
\eeq

\noi for $0 \leq \ell \leq k-1$. Let now $0 < t_1 < t_2 \leq \tau_0$. Integrating (\ref{4.44e}) over time yields
\beq
\label{4.45e}
\parallel \widetilde{v} (t_2) - \widetilde{v} (t_1) ; H^{k-1} \parallel \ \leq C\ Y^3 (1 + Y^2) t_2^{k/2} \ .
\eeq

\noi From (\ref{4.45e}) it follows that $\widetilde{v}  (t)$ has a limit $v_+$ in $H^{k-1}$ when $t\to 0$ and converges to that limit according to  (\ref{4.34e}). Together with uniform boundedness of $v$ in $V^k$, this implies that $v_+ \in V^k$, that $v_+$ satisfies  (\ref{4.33e}), and that $v(t)$ converges to $v_+$ in the convergences stated in Part (1). Similar arguments apply with $V^k$ replaced by $\Sigma^k$.\\

\noi \underline{Part (2).} From  (\ref{2.21e})  (\ref{4.37e}) we obtain 
\beq
\label{4.46e}
\partial_t (\varphi - \varphi_1) = {\check B}_L (v - v_a, v + v_a) \ .
\eeq

\noi Using (\ref{3.43e}) and Lemma 3.3, we estimate 
$$\parallel   \omega^{\ell} \partial_t (\varphi - \varphi_1)\parallel _2\ \leq\ C \ t^{-(\ell / 2 + 1 - k)_+}\parallel  \omega^{2(k-1) \wedge \ell} \partial_t (\varphi - \varphi_1) \parallel _2 $$ 
$$\leq C\ t^{-1-(\ell /2 + 1 - k)_+}  \parallel v - v_a ; H^{k-1} \parallel \ \parallel v + v_a ; V^k \parallel$$
\beq
\label{4.47e}
\leq C\ Y^4 (1 + Y^2) t^{-1 + k/2 - (\ell /2 + 1 - k)_+} 
\eeq

\noi for $\ell > 0$. Integrating (\ref{4.47e}) over time yields
\beq
\label{4.48e}
\parallel   \omega^{\ell} \left ( \varphi (t_2) - \varphi_1 (t_2 ) - \varphi (t_1) + \varphi_1(t_1)\right ) \parallel_2\ \leq \ C\ Y^4 (1 + Y^2) t_2^{k/2-(\ell /2 + 1 - k)_+}
\eeq

\noi for $0 < \ell < 2k-3$, which implies the existence of $\psi_{1+}$ and the estimate (\ref{4.38e}). The estimate (\ref{4.39e}) follows from (\ref{4.18e}) (\ref{4.38e}), from the inequality
$$k/2 - (\ell /2 + 1 - k)_+ \geq (k-1)/2 - (\ell /2 + 1/2 - k)_+$$

\noi and from the fact that $Y^2 t^{(k-1)/2} \leq C$ by (\ref{4.17e}) (\ref{4.19e}).\\

\noi \underline{Part (3).} From (\ref{4.37e}) (\ref{4.40e}) it follows that
\beq
\label{4.49e}
\partial_t \left ( \varphi_1 - \varphi_2 \right ) = {\check B}_L \left ( v_a - v_+, v_a + v_+\right ) - {\check B}_S(v_+)  \ .
\eeq

\noi We estimate
\bea
\label{4.50e}
&&\parallel   \omega^{\ell} {\check B}_L \left ( v_a - v_+, v_a + v_+\right )\parallel_2\ \leq \ C\ t^{-1} \parallel  v_a - v_+ \parallel _2 \ \parallel  v_a+v_+; V^k \parallel   \nn \\
&&\qquad \leq C\ Y^2\ t^{-1+k/2} \qquad \hbox{for $0 < \ell \leq k-1$}
\eea

\noi by Lemma 3.3, which together with (\ref{3.43e}) yields 
\beq
\label{4.51e}
\parallel   \omega^{\ell} {\check B}_L \left ( v_a - v_+, v_a + v_+\right )\parallel_2\ \leq \ C\ Y^2\ t^{-1+k/2-(1/2) (\ell + 1 -k)_+}
\eeq

\noi for $\ell > 0$. On the other hand
\beq
\label{4.52e}
\parallel   \omega^{\ell} {\check B}_S (v_+)\parallel_2\ \leq \ C\ Y^2\ t^{-1+k- \ell /2 - 1/2}\qquad \hbox{for $\ell \leq 2k-1$}
\eeq

\noi by Lemma 3.3 and (\ref{3.42e}), so that 
\beq
\label{4.53e}
\parallel   \omega^{\ell} \partial_t \left ( \varphi_1 - \varphi_2 \right )\parallel_2\ \leq \ C\ Y^2\ t^{-1+k/2-(1/2) (\ell + 1 -k)_+}
\eeq

\noi for $0 < \ell \leq 2k-1$. Combining  (\ref{4.53e}) with  (\ref{4.47e}) and using the same arguments as in the proof of Part (2) yields the existence of $\psi_{2+}$ and the estimate (\ref{4.41e}), which together with (\ref{4.18e}) yields (\ref{4.42e}). \par \nobreak \hfill $\sq$\par

\noi {\bf Remark 4.2.} From (\ref{4.44e}) one can obtain slightly better estimates of  $\parallel   \omega^{\ell} (v-v_a)\parallel_2$ for $0 \leq \ell < k-1$, which imply slightly better estimates of  $\parallel   \omega^{\ell} (\varphi -\varphi_a)\parallel_2$ for $0 < \ell < 2(k-1)$.\\

We now turn to the Cauchy problem with finite initial time for the original system (\ref{1.1e}) and we state the results on that problem that follow from the results for the auxiliary system contained in Propositions 4.2 and 4.3. We recall that $\widetilde{u}(t) = U(-t) u(t) = \overline{F \widetilde{u}_c(1/t)}$.\\

\noi {\bf Proposition 4.4.} {\it Let $1 < k < 2$ and let $\widetilde{u}_0 \in FV^k$ with $\parallel \widetilde{u}_0;FV^k \parallel \ = \widetilde{a}$. Then\par

(1) There exists $\overline{t}_0 \geq 1$ such that for any $t_0 \geq  \overline{t}_0$ there exists a unique solution $u$ of the system (\ref{1.1e}) with $\widetilde{u} \in {\cal C} (I, FV^k)$ and $u(t_0) = U(t_0)\widetilde{u}_0$, where $I = [t_0, \infty )$. The time $\overline{t}_0$ depends on $\widetilde{a}$ according to
\beq
\label{4.54e}
\widetilde{a} \leq C\ \overline{t}_0^{(k-1)/4}\ .
\eeq

\noi In particular one can take $\overline{t}_0 = 1$ for small $\widetilde{a}$. The solution $u$ satisfies the following estimate
\beq
\label{4.55e}
\parallel \widetilde{u}(t); FV^k \parallel \ \leq C\ \widetilde{a} \left ( 1 + \widetilde{a}^2(1 + \ell n\ t)\right )^{k \vee 3/2}   
\eeq

\noi for all $t\in I$. Define in addition $\varphi$ and $\theta$ by
\beq
\label{4.56e}
\varphi (1/t_0) = 0 \quad , \quad \partial_t \varphi = {\check B}_L(u_c) \qquad \hbox{for $0 < t \leq 1/t_0$}\ ,
\eeq
\beq
\label{4.57e}
\theta (t) = - D_0 (t) \varphi (1/t) \qquad \hbox{for $t \geq t_0$}\ .
\eeq

\noi Then the following estimates hold for all $t \in I$
\beq
\label{4.58e}
\parallel   \omega^{\ell} \theta (t)\parallel_2\ \leq \ C\ \widetilde{a}^2 \ t^{1- \ell } \left ( t^{(\ell /2 + 1/2 -k)_+} + \ell n \ t\right ) \qquad \hbox{for $\ell > 0$}\ ,
\eeq
\beq
\label{4.59e}
\parallel  U(-t) u(t) \exp (i \theta (t));L^{\infty} (I, FV^k)\parallel  \ \leq \ C\ \widetilde{a}\ .
\eeq

\noi If in addition $\widetilde{u}_0 \in \Sigma^k$, then $u$, $\widetilde{u} \in {\cal C}(I, \Sigma^k)$ and $u$ satisfies the estimates
\beq
\label{4.60e}
\parallel \widetilde{u}(t); \Sigma^k \parallel \ \leq C\ \widetilde{a} \left ( 1 + \widetilde{a}^2(1 + \ell n(1+ |t|))\right )^{k \vee 3/2}   
\eeq

\noi for all $t \in I$, 
\beq
\label{4.61e}
\parallel  U(-t) u(t) \exp (i \theta (t));L^{\infty} (I, \Sigma^k)\parallel  \ \leq \ C\ \widetilde{a}
\eeq

\noi where now $\widetilde{a} = \ \parallel \widetilde{u}_0; \Sigma^k\parallel $.\par

(2) The map $\widetilde{u}_0 \to \widetilde{u}$ is continuous for fixed $t_0$ on the bounded sets of $FV^k$ from the $L^2$ norm of $\widetilde{u}_0$ to the norm of $\widetilde{u}$ in $L^{\infty} (J, FH^{k'})$ for $0 \leq k' < k$ and in the weak $\star$  sense in $L^{\infty} (J, FV^{k})$ for any interval $J \subset\subset I$. If in addition $\widetilde{u}_0 \in \Sigma^k$, then continuity holds on the bounded sets of $\Sigma^k$ to the norm of $\widetilde{u}$ in $L^{\infty} (J, \Sigma^{k'})$ for $0 \leq k' < k$ and in the weak $\star$ sense in $L^{\infty} (J, \Sigma^{k})$. \par

(3) Let $u_0 \in \Sigma^k$ with $\widetilde{a} = \ \parallel u_0; \Sigma^k \parallel$ sufficiently small. Then there exists a unique solution $u$ of the system (\ref{1.1e}) with $u$, $\widetilde{u} \in {\cal C} ({I\hskip-1truemm R}, \Sigma^k)$ and $u(0) = u_0$. That solution satisfies (\ref{4.60e}) for all $t \in {I\hskip-1truemm R}$.}\\

\noi {\bf Proof.} \underline{Part (1).} We first prove the existence of a solution with the properties stated. Let $\tau_0 = 1/t_0$ and 
\beq
\label{4.62e}
v_0 = U(\tau_0) \widetilde{v}_0 = U(1/t_0) \overline{F\widetilde{u}_0}\ .
\eeq

\noi Let $(v, \varphi )$ be the solution of the system (\ref{2.21e}) (\ref{2.22e}) obtained in Proposition 4.2 with $(v, \varphi ) (\tau_0) = (v_0, 0)$. Such a solution exists for $\tau_0 \leq \overline{\tau}_0$ and $\overline{\tau}_0$ satisfying (\ref{4.19e}). Now by (\ref{3.8e})
\beq
\label{4.63e}
a = \ \parallel v_0;V^k\parallel\ \leq\ 2 \parallel \widetilde{v}_0; V^k \parallel\ = 2 \ \parallel \widetilde{u}_0; FV^k \parallel\ = 2 \widetilde{a}
\eeq

\noi so that (\ref{4.19e}) follows from (\ref{4.54e}) (with a different constant). Define $u$ by  (\ref{2.6e}) (\ref{2.16e}). Then $u$ solves  (\ref{1.1e}) in $I = [t_0, \infty )$ with $u(t_0) = U(t_0)  \widetilde{u}_0$, and $\varphi$ satisfies (\ref{4.56e}) because ${\check B}_L (v) = {\check B}_L (u_c)$, so that $\varphi$ can actually be defined in terms of $u$ by (\ref{4.56e}). Furthermore
\beq
\label{4.64e}
U(-t) \ u(t) \exp (i \theta (t)) = \overline{F\widetilde{v}(1/t)}\ .
\eeq   

\noi The regularity of $u$ follows immediately from that of $(v, \varphi )$ through (\ref{2.6e}) (\ref{2.16e}). The estimates (\ref{4.58e}) (\ref{4.59e}) are essentially a rewriting of (\ref{4.18e})  (\ref{4.17e}). We next derive (\ref{4.55e}). Now by (\ref{3.8e}) 
\bea
\label{4.65e}
\parallel \widetilde{u}(t); FV^k \parallel &=& \parallel  \widetilde{u}_c (1/t); V^k\parallel \ \leq \ 2\parallel u_c (1/t); V^k \parallel\nn \\
 &=&2 \parallel \left ( v \exp (- i \varphi )\right )(1/t);V^k \parallel  
\eea

\noi and we estimate the last norm by using (\ref{4.17e}) (\ref{4.18e}) (\ref{4.63e}) and Lemma 3.4, part (1) with $m = (2k-1) \wedge 2$. This proves (\ref{4.55e}).\par

We finally prove uniqueness of $u$ by estimating the $L^2$ norm of the difference of the pseudo conformal inverses $u_{ci}$, $i=1,2$, of two solutions $u_i$, $i=1,2$. From (\ref{2.9e}), by a simplified version of Lemma 4.2, part (1), we estimate
\beq
\label{4.66e}
\left | \partial_t \parallel u_{c-} \parallel_2 \right | \ \leq C\left ( t^{-1}y^2 + y^4\right )  \parallel u_{c-}\parallel_2
\eeq

\noi where
$$y = y(t) = \ \mathrel{\mathop {\rm Max}_i}\ \parallel u_{ci} (t); V^k \parallel$$

\noi from which uniqueness follows immediately.\par

The additional properties of $u$ for $\widetilde{u}_0 \in \Sigma^k$ follow immediately from the last statement of Proposition 4.2, part (1) by similar arguments.\\

\noi \underline{Part (2)} follows immediately from Proposition 4.2, part (2).\\

\noi \underline{Part (3).} By (\ref{4.54e}), for $ \widetilde{a}$ sufficiently small, we can take $\overline{t}_0 = 1$ in Part (1) of this proposition. Applying that result, we obtain a solution $u$ of the system (\ref{1.1e}) with $u(1) = u_0$ and $ \widetilde{u} \in {\cal C}([1, \infty ) , \Sigma^k)$ provided
\beq
\label{4.67e}
 \widetilde{a}_> = \ \parallel U(-1) u_0 ; \Sigma^k \parallel \ \leq \overline{a}
\eeq

\noi for some $\overline{a}$ sufficiently small. Since the system (\ref{1.1e}) is time translation invariant, by translating the previous solution by $-1$ in time, we obtain a solution $u_>$ with $u_> (0) = u_0$ and $U(-1)  \widetilde{u}_> \in {\cal C} ([0, \infty ) , \Sigma^k)$, or equivalently $\widetilde{u}_> \in {\cal C} ([0, \infty ) , \Sigma^k)$, satisfying the estimate (\ref{4.60e}) for all $t \geq 0$ with $ \widetilde{a}$ replaced by $ \widetilde{a}_>$. Since the system (\ref{1.1e}) is also time reversal invariant, we can construct similarly a solution $u_<$ with $u_<(0) = u_0$ and $\widetilde{u}_< \in {\cal C} (( - \infty , 0]  , \Sigma^k)$, satisfying the estimate (\ref{4.60e}) for all $t \leq 0$ with $ \widetilde{a}$ replaced by $ \widetilde{a}_<$, with 
\beq
\label{4.68e}
\widetilde{a}_< = \ \parallel U(1) u_0 ; \Sigma^k \parallel \ \leq \overline{a} \ .
\eeq

\noi Taking $u(t) = u {>\atop <}(t)$ for $t{>\atop <} 0$ yields a solution u of the system  (\ref{1.1e}) with $u(0) = u_0$ and $\widetilde{u} \in {\cal C} ({I\hskip-1truemm R}, \Sigma^k)$, satisfying the estimate (\ref{4.60e}) for all $t \in {I\hskip-1truemm R}$ with $\widetilde{a} = \widetilde{a}_> \vee \widetilde{a}_<$. Finally by (\ref{3.11e}), the conditions (\ref{4.67e}) (\ref{4.68e}) can both be satisfied by taking \break\noindent $\parallel u_0; \Sigma^k \parallel$ sufficiently small.\par \nobreak \hfill $\sq$ \par

\noi {\bf Remark 4.3.} Proposition 4.3 allows one to take arbitrarily large $\widetilde{u}_0 \in \Sigma^k$ by taking $\overline{t}_0$ sufficiently large according to (\ref{4.54e}), thereby generating some large initial data $u_0 = U(t_0) \widetilde{u}_0$ in $\Sigma^k$. However one cannot accomodate arbitrarily large $u_0 \in \Sigma^k$ since for fixed $u_0 \in \Sigma^k$ and $t_0$ large 
$$\parallel <x>^k \widetilde{u}_0\parallel_2 \ = \ \parallel <x+it_0 \nabla >^k u_0 \parallel_2 \ \sim \ t_0^k \parallel u_0 ; H^k \parallel$$

\noi and taking $t_0$ large is of no help in order to fulfill (\ref{4.54e}).\\

We now turn to the asymptotic properties of the solutions obtained in Proposition 4.4 that follow from Proposition 4.3.\\

\noi {\bf Proposition 4.5.} {\it Let $1 < k <2$. Let $u$ be a solution of the system (\ref{1.1e}) as obtained in Proposition 4.4, let $\theta$ be defined by (\ref{4.56e}) (\ref{4.57e}) and let 
\beq
\label{4.69e}
\widetilde{Y}=\  \parallel U(-t) \ u(t) \exp ( i \theta (t)) ;L^{\infty}(I, FV^k  \parallel
\eeq

\noi where $I = [t_0 , \infty )$. Then \par

(1) There exists $u_+ \in FV^k$ such that $U(-t) u(t)  \exp ( i \theta (t))$ tends to $u_+$ when $t \to \infty$ strongly in $FV^{k, \rho}$ for $0 \leq \rho < 1$ (and in particular in $FH^{k-1}$) and weakly in $FV^k$. Furthermore 
\beq
\label{4.70e}
\parallel u_+; FV^k \parallel \ \leq \ \mathrel{\mathop {\rm lim\ inf}_{t\to \infty}}\  \parallel U(-t)\ u(t) \exp ( i \theta (t)) ;FV^k  \parallel\ \leq \widetilde{Y}
\eeq

\noi and the following estimate holds for all $t \in I$~:
\beq
\label{4.71e}
\parallel <x>^{k-1} \left ( U(-t)\ u(t) \exp (i \theta (t)) - u_+ \right )  \parallel_2 \ \leq C\   \widetilde{Y}^3 (1 + \widetilde{Y}^2) t^{-k/2}\ .
\eeq

\noi If in addition $\widetilde{u} \in {\cal C}(I, \Sigma^k)$, then $u_+ \in \Sigma^k$ and $U(-t) u(t) \exp (i \theta (t))$ tends to $u_+$ when $t \to \infty$ strongly in $\Sigma^{k'}$ for $0 \leq k'< k$ and weakly in $\Sigma^k$. Furthermore
\beq
\label{4.72e}
\parallel u_+; \Sigma^k \parallel \ \leq \ \mathrel{\mathop {\rm lim\ inf}_{t\to \infty}}\  \parallel U(-t)\ u(t) \exp ( i \theta (t)) ; \Sigma^k  \parallel \ .
\eeq

(2) Let $\varphi_a$ be defined as in Proposition 4.3, part (2) with $v_+ = \overline{Fu_+}$ and define
\bea
\label{4.73e}
&&\theta_a (t) = - D_0 (t) \ \varphi_a (1/t)\ , \\
&&u_a(t) = \exp (-i \theta_a (t)) U(t)\ u_+ \ .
\label{4.74e}
\eea

\noi Then $u_a$ satisfies the estimate
\beq
\label{4.75e}
\parallel \widetilde{u}_a(t); F V^k \parallel \ \leq C\ \widetilde{Y}\left ( 1 + (\widetilde{Y}^2(1 + \widetilde{Y}^2)(1 + \ell n\ t))^{k\vee 3/2}\right ) \ .
\eeq

\noi Furthermore $u$ behaves asymptotically as $u_a$ for large $t$ in the sense that $\widetilde{u} - \widetilde{u}_a$ tends to zero when $t \to \infty$ strongly in $FV^{k, \rho}$ for $0 \leq \rho < 1$ (and in particular in $FH^{k-1})$. The difference $\widetilde{u} - \widetilde{u}_a$ satisfies the estimate 
\beq
\label{4.76e}
\parallel <x>^{k-1} \left ( \widetilde{u}(t) -  \widetilde{u}_a(t)\right ) \parallel_2  \ \leq C\ \widetilde{Y}^3( 1 + \widetilde{Y}^2)^3(1 + \ell n\ t) t^{-k/2}
\eeq

\noi for all $t \in I$, and similar estimates in $FV^{k , \rho}$, $0 \leq \rho < 1$, obtained by interpolation between (\ref{4.76e}) and the estimates (\ref{4.55e}) (\ref{4.75e}) in $FV^k$. If in addition $\widetilde{u} \in {\cal C}(I, \Sigma^k)$, then $\widetilde{u}_a$ satisfies an estimate in $\Sigma^k$ similar to (\ref{4.75e}) and $\widetilde{u} - \widetilde{u}_a$ tends to zero when $t\to \infty$ strongly in $\Sigma^{k'}$ for $0 \leq k' < k$.\par

(3) Let $\varphi_b$ be defined as in Proposition 4.3, part (3) with $v_+ = \overline{Fu_+}$ and define
\beq
\label{4.77e}
\theta_b (t) = - D_0 (t) \ \varphi_b (1/t) = D_0 (t) \left ( (\ell n\ t ) x \cdot B (\overline{Fu_+}) - \psi_{2+}\right )\ , 
\eeq
\beq
\label{4.78e}
u_b (t) = \exp \left (- i \theta_b (t)\right ) U(t)\ u_+ \ .
\eeq

\noi Then $u_b$ satisfies the estimate
\beq
\label{4.79e}
\parallel \widetilde{u}_b(t); F V^k \parallel \ \leq C_{\varepsilon} \ \widetilde{Y}\left ( 1 + \left ( \widetilde{Y}^2 (1 + \widetilde{Y}^2)^2 (1 + \ell n\ t)\right ) ^{k\vee (3/2+ \varepsilon )}\right ) 
\eeq

\noi for any $\varepsilon > 0$. Furthermore $u$ behaves asymptotically as $u_b$ for large $t$ in the sense that $\widetilde{u} - \widetilde{u}_b$ tends to zero when $t\to \infty$ in $FV^{k, \rho}$, $0 \leq \rho < 1$ (and in particular in $FH^{k-1}$). The difference $\widetilde{u} - \widetilde{u}_b$ satisfies the estimate
\beq
\label{4.80e}
\parallel <x>^{k-1} \left ( \widetilde{u}(t) -  \widetilde{u}_b(t)\right ) \parallel_2  \ \leq C\ \widetilde{Y}^3( 1 + \widetilde{Y}^2)^3(1 + \ell n\ t) t^{-k/2} \ .
\eeq

\noi for all $t \in I$, and similar estimates in $FV^{k, \rho }$, $0 \leq \rho < 1$, obtained by interpolation between (\ref{4.80e}) and the estimates (\ref{4.55e}) (\ref{4.79e}) in $FV^k$.\par

If $\widetilde{u} \in {\cal C} (I, \Sigma^k)$ a similar reinforcement occurs as in Part (2).}\\

\noi {\bf Proof.} \underline{Part (1).} The solution $u$ is obtained from a solution $(v, \varphi )$ of the system (\ref{2.21e}) (\ref{2.22e}) as in the proof of Proposition 4.4 and $u$ satisfies (\ref{4.59e}) with $\theta$ defined by (\ref{4.56e}) (\ref{4.57e}).\par

The existence of the limit $u_+$ and the convergence properties of Part (1) are a rewriting of the corresponding properties in Proposition 4.3, part (1) with $u_+ = \overline{Fv_+}$. The estimates (\ref{4.70e}) (\ref{4.71e}) (\ref{4.72e}) follow from (\ref{4.33e}) (\ref{4.34e})  (\ref{4.36e}), from the relations
\beq
\label{4.81e}
\parallel U(-t) \ u(t)  \exp (i \theta (t)) - u_+ ; FH^{k-1} \parallel \ =\   \parallel  \widetilde{v}(1/t) - v_+ ; H^{k-1} \parallel \ , 
\eeq
\beq
\label{4.82e}
Y\ = \ \parallel v;L^{\infty} ((0, \tau_0];V^k) \parallel \ \leq \   2\parallel  \widetilde{v}; L^{\infty} ((0, \tau_0], V^k)\parallel \ = \ 2 \widetilde{Y} \ , 
\eeq
  
\noi by (\ref{3.8e}) (\ref{4.32e}) (\ref{4.69e}), and a similar one with $V^k$ replaced by $\Sigma^k$. \\

\noi \underline{Part (2).} The definition (\ref{4.74e}) of $u_a$ is actually a rewriting of 
$$\widetilde{u}_a (t) = \overline{F \widetilde{u}_{ca} (1/t)}$$

\noi with
$$u_{ca} = v_a \exp (- i \varphi_a) \qquad , \quad v_a (t) = U(t) v_+ \ ,$$

\noi in analogy with  (\ref{2.6e}) (\ref{2.16e}). Therefore
\beq
\label{4.83e}
\parallel \widetilde{u}_a (t) ; FV^k \parallel \ \leq \   2\parallel  \left ( v_a \exp (- i \varphi_a ) \right ) (1/t); V^k \parallel 
\eeq

\noi by (\ref{3.8e}), which implies (\ref{4.75e}) by Lemma 3.4, part (1) with $m = (2k-1)\wedge 2$, (\ref{4.33e}) (\ref{4.39e}) (\ref{4.82e}). Similarly
\beq
\label{4.84e}
\parallel \widetilde{u} (t)- \widetilde{u}_a (t) ; FH^{k-1} \parallel \ =  \   \parallel  \left ( v \exp (- i \varphi ) - v_a (\exp - i \varphi_a)\right ) (1/t); H^{k-1} \parallel \ . 
\eeq
\noi In order to prove (\ref{4.76e}), we write
\begin{eqnarray*}
&&v \exp (- i \varphi ) - v_a \exp (- i \varphi_a) = v_{\not=} \exp (- i \varphi )\ ,\\
&&v_{\not=}  = v - v_a + v_a (1 - \exp (i \psi ))
\end{eqnarray*} 

\noi with $\psi = \varphi - \varphi_a$ and for $0 < \ell \leq k-1$, we estimate
\beq
\label{4.85e}
\parallel   \omega^{\ell} v_{\not=} \exp (-i \varphi ) \parallel_2\ \leq \ C \parallel   \omega^{\ell} v_{\not=}  \parallel_2 \left ( 1 + \ \parallel \nabla \varphi \parallel_2\right )\ ,
\eeq
\beq
\label{4.86e}
\parallel   \omega^{\ell} v_{\not=} \parallel_2\ \leq \ C  \left (\parallel   \omega^{\ell} (v- v_a)  \parallel_2 \  + \left (  \parallel v_a \parallel_{\infty} \ + \  \parallel  \nabla v_a  \parallel _2 \right )  \parallel \omega^{\ell}  \psi \parallel_2\right ) 
\eeq

\noi by Lemma 3.1 and (\ref{3.16e}), which implies that $\parallel \omega^{\ell} \exp (i \psi )\parallel_2 \ \leq \ \parallel \omega^{\ell}\psi \parallel_2$. We continue (\ref{4.85e}) (\ref{4.86e}) by using (\ref{4.17e}) (\ref{4.18e}) (\ref{4.33e}) (\ref{4.34e}) (\ref{4.38e}) (\ref{4.82e}). Substituting the result into (\ref{4.84e}) yields (\ref{4.76e}), which by interpolation with (\ref{4.55e}) (\ref{4.75e}) completes the proof in the case of $FV^k$.\par

The proof in the case of $\Sigma^k$ is similar.\\

\noi \underline{Part (3)}. The proof of (\ref{4.79e}) is essentially the same as that of (\ref{4.75e}) with the estimate (\ref{4.39e}) replaced by (\ref{4.42e}). The occurrence of $\varepsilon > 0$ in the exponent for $k \leq 3/2$ is required by the fact that (\ref{4.42e}) holds only for $\ell < 2k-1$, the limiting case being excluded. Similarly the proof of (\ref{4.80e}) is essentially the same as that of (\ref{4.76e}), with (\ref{4.38e}) replaced by (\ref{4.41e}). Note in particular that (\ref{4.38e}) (\ref{4.41e}) are used in those proofs with $\ell \leq k-1$ and yield the same time decay in that case. \par \nobreak \hfill $\sq$ \par

\noi {\bf Remark 4.4.} In Parts 2 and 3 of Proposition 4.5, the convergence of $\widetilde{u} - \widetilde{u}_a$ to zero as $t \to \infty$ can be extended to weak convergence in $FH^k$ for solutions in  ${\cal C}(I, FV^k)$ and to weak convergence in $\Sigma^k$ for solutions in ${\cal C}(I, \Sigma^k)$. This follows from convergence in $FH^{k-1}$ and from uniform boundedness of $\widetilde{u} - \widetilde{u}_a$ in $FH^k$ or in $\Sigma^k$. The latter follows from uniform boundedness of $v_{\not=}\exp (- i \varphi )$ in $H^k$ or in $\Sigma^k$, which can be easily derived from the available estimates. 

\mysection{The Cauchy problem at initial time zero for\break\noindent (v,$\varphi$) and at infinity for u}
\hspace*{\parindent} In this section we solve the Cauchy problem for the auxiliary system (\ref{2.21e}) (\ref{2.22e}) with prescribed asymptotic behaviour at time zero and for the original system (\ref{1.1e}) with corresponding prescribed asymptotic behaviour at infinity. The main results are contained in Proposition 5.6 for the system (\ref{2.21e}) (\ref{2.22e}) and in Proposition 5.7 for the system (\ref{1.1e}). In this section we use only the spaces $V^k$ for $v$ and $FV^k$ for $\widetilde{u}$. The specialization of the results to the space $\Sigma^k$ is straightforward and will not be considered. We start with a uniqueness result for $(v, \varphi )$.\\

\noi {\bf Proposition 5.1.} {\it Let $1 < k <2$. Let $0 < \tau \leq 1$ and $I = (0, \tau ]$. Let $(v_i, \varphi_i)$, $i = 1,2$, be two solutions of the system (\ref{2.21e}) (\ref{2.22e}) with $v_i \in {\cal C}(I, V^k)$ and $\varphi_i (t_0) \in H_>^{k+2}$ for some $t_0 \in I$. Assume that 
\beq
\label{5.1e}
\parallel v_i(t) ; V^k \parallel \ \leq a (1 - \ell n\ t)^{\alpha}
\eeq

\noi for all $t \in I$ and 
\beq
\label{5.2e}
\mathrel{\mathop {\rm Sup}_{t\in I}} \ h_1(t)^{-1} \parallel v_1(t) - v_2(t) \parallel_2\ = Y < \infty 
\eeq

\noi for some constants $a$, $Y$ and $\alpha \geq 0$ and for some non decreasing function $h_1 \in {\cal C} (I, {I\hskip-1truemm R}^+)$ satisfying
\beq
\label{5.3e}
\int_0^t dt'\ t{'}^{-1 -(3-k)/2} (1 - \ell n\ t ')^{\alpha} \ h_1(t') \leq C\ t^{-(3-k)/2} (1 - \ell n \ t)^{\alpha} \ h_1(t)
\eeq

\noi for all $t \in I$, and
\beq
\label{5.4e}
\lim_{t \to 0} t^{-(3-k)/2} (1 - \ell n\ t )^{\alpha} \ h_1(t) = 0 \ .
\eeq

\noi Assume in addition that 
\beq
\label{5.5e}
\lim_{t \to 0} \left ( \varphi_1(t) - \varphi_2 (t) \right ) = 0\ .
\eeq

\noi Then $(v_1, \varphi_1) = (v_2, \varphi_2 )$.}\\

\noi {\bf Proof.} We define again $(v_{\pm} , \varphi_{\pm}) = (1/2) (v_1 \pm v_2, \varphi_1 \pm \varphi_2)$ and we estimate $(v_-, \varphi_-)$ by Lemma 4.2. We first estimate $\varphi_-$. From (\ref{4.15e}) (\ref{5.1e}) (\ref{5.2e}) it follows that 
\beq
\label{5.6e}
\parallel \omega^{\ell} \partial_t \varphi_- \parallel_2 \ \leq \ C\ a \ Y(1 - \ell n\ t)^{\alpha} \ t^{-1-(1/2)(\ell + 1 - k)_+} \ h_1(t)
\eeq

\noi for $\ell > 0$. From (\ref{5.3e}) (\ref{5.4e}) (\ref{5.6e}) it follows that $\varphi_-$ has a limit in $H_>^2$ when $t \to 0$, which gives a meaning to the assumption (\ref{5.5e}). Furthermore
\beq
\label{5.7e}
\parallel \omega^{\ell}  \varphi_-(t) \parallel_2 \ \leq \ C\ a \ Y(1 - \ell n\ t)^{\alpha} \ t^{-(1/2)(\ell + 1 - k)_+} \ h_1(t)
\eeq 

\noi for $0 < \ell \leq 2$ and all $t \in I$. On the other hand, by integrating (\ref{4.6e}) between $t_0$ and $t$, we obtain 
\beq
\label{5.8e}
\parallel \omega^{\ell}  \varphi_i(t) \parallel_2 \ \leq \ C\ a^2 (1 - \ell n\ t)^{2\alpha} \left ( t^{-(\ell /2 + 1/2 - k)_+} - \ell n\ t \right )
\eeq 

\noi for $\ell > 0$. We next estimate $v_-$ by Lemma 4.2. From (\ref{4.14e}) (\ref{5.1e}) (\ref{5.8e}) we obtain
$$ \left | \partial_t \parallel v_{-} \parallel_2 \right | \ \leq C\  a (1 - \ell n\ t)^{\alpha} \Big \{ \parallel \omega^{3-k} \varphi_-\parallel_2\ +\ \parallel \omega^2 \varphi_-\parallel_2$$
$$+ \ a (1 - \ell n\ t)^{\alpha}\ t^{-(3-k)/2} \parallel  v_- \parallel _2 \ + \ a^2 (1 - \ell n\ t)^{2\alpha}\left ( t^{-(3/2  - k)_+} - \ell n\ t \right )\parallel \omega \varphi_-\parallel _2$$
\beq
\label{5.9e}
+ \ a^3 (1 - \ell n\ t )^{3\alpha + 1} \parallel v_- \parallel_2 \Big \} \ .
\eeq

\noi From (\ref{5.9e}) (\ref{5.2e}) (\ref{5.7e}) we then obtain
\bea
\label{5.10e}
 \left | \partial_t \parallel v_{-} \parallel_2 \right | &\leq& C\  a^2 (1 - \ell n\ t)^{2\alpha} \  t^{-(3 - k)/2} \ Y \ h_1(t) \nn\\
&\times&  \left ( 1 + a^2(1 - \ell n\ t)^{2\alpha}\ t^{1/2} \left ( t^{-(3/2  - k)_+} - \ell n\ t \right )\right ) \nn \\
&\leq& C\  a^2 (1 + a^2) (1 - \ell n\ t )^{2\alpha} \ t^{-(3 - k)/2} \ Y \ h_1(t) \ .
\eea

\noi Integrating (\ref{5.10e}) in $(0, t]$ and using (\ref{5.2e}) (\ref{5.3e}), we obtain
\beq
\label{5.11e}
Y \leq C\ a^2 (1 + a^2) (1 - \ell n\ \tau)^{2\alpha} \ \tau^{(k-1)/2}\ Y
\eeq

\noi which implies $Y = 0$ by taking $\tau$ sufficiently small and therefore $v_- = 0$ and $\varphi_-=0$ by (\ref{5.7e}). The extension of the proof to larger $\tau$ proceeds by standard arguments.\par\nobreak\hfill $\sq$ \par

We next derive a uniqueness result for solutions of the system (\ref{1.1e}). That result is obtained by applying Proposition 5.1 to solutions of the system (\ref{2.21e})(\ref{2.22e}) reconstructed from solutions of the system (\ref{1.1e}).\\

\noi {\bf Proposition 5.2.} {\it Let $1 < k <2$. Let $T \geq 1$ and $I = [T , \infty )$. Let $u_i$, $i = 1,2$ be two solutions of the system (\ref{1.1e}) with $\widetilde{u}_i \in {\cal C}(I, FV^k)$, where $\widetilde{u}(t) = U(-t) u(t)$. Assume that
\beq
\label{5.12e}
\parallel \widetilde{u}_i (t); FV^k \parallel\ \leq b (1 + \ell n\ t)^{\beta}
\eeq

\noi for all $t \in I$ and
\beq
\label{5.13e}
\mathrel{\mathop {\rm Sup}_{t\in I}} \ h_2(1/t)^{-1} \parallel \widetilde{u}_1(t) - \widetilde{u}_2(t) \parallel_2\ = Z < \infty 
\eeq

\noi for some constants $b$, $Z$ and $\beta \geq 0$ and for some nondecreasing function $h_2 \in {\cal C}(I_0 , {I\hskip-1truemm R}^+)$ with $I_0 = (0, T^{-1}]$, satisfying
\beq
\label{5.14e}
\int_0^t dt'\ t{'}^{-1 -(3-k)/2} (1 - \ell n\ t ')^{7 \beta + 2} \ h_2(t') \leq C\ t^{-(3-k)/2} (1 - \ell n \ t)^{7 \beta + 2} \ h_2(t)
\eeq

\noi for all $t\in I_0$, and
\beq
\label{5.15e}
\lim_{t \to 0} t^{-(3-k)/2} (1 - \ell n\ t )^{7\beta +2} h_2(t) = 0 \ .
\eeq

\noi Then $u_1 = u_2$.} \\

\noi {\bf Proof.} From $u_i$, $i = 1,2$, we reconstruct solutions $(v_i, \varphi_i)$ of the system (\ref{2.21e}) (\ref{2.22e}). For convenience, we work with the pseudoconformal inverses $u_{ci}$, $i = 1, 2$ of $u_i$, defined by (\ref{2.6e}). From (\ref{3.8e}) it follows that (\ref{5.12e}) (\ref{5.13e}) can be rewritten as
\beq
\label{5.16e}
\parallel u_{ci} (t); V^k \parallel\ \leq b (1 - \ell n\ t)^{\beta}
\eeq
\beq
\label{5.17e}
\mathrel{\mathop {\rm Sup}_{t\in I_0}} \ h_2(t)^{-1} \parallel  u_{c1}(t) - u_{c2}(t) \parallel_2\ = Z < \infty 
\eeq

\noi possibly with a change of $b$ by a factor 2. We first define the phases $\varphi_i$. We define $\varphi_2$ by 
\beq
\label{5.18e}
\partial_t \ \varphi_2 = {\check B}_L (u_{c2})
\eeq

\noi with initial condition $\varphi_2 (t_0) = 0$ for some $t_0 \in I_0$. By integrating (\ref{4.6e}) in $[t_0 , t]$ and using (\ref{5.16e}) we obtain
\beq
\label{5.19e}
\parallel \omega^{\ell}  \varphi_2(t) \parallel_2 \ \leq \ C\ b^2 (1 - \ell n\ t)^{2\beta} \left ( t^{-(\ell /2 + 1/2 - k)_+} - \ell n\ t \right )
\eeq 

\noi for $\ell > 0$ and in particular
\beq
\label{5.20e}
\parallel \omega^{\ell}  \varphi_2(t) \parallel_2 \ \leq \ C\ b^2 (1 - \ell n\ t)^{2\beta + 1} \qquad \hbox{for $0 < \ell \leq 2k-1$}\ .
\eeq  

\noi We next define $\varphi_-$ as a solution of the equation 
\beq
\label{5.21e}
\partial_t \ \varphi_- = {\check B}_L (u_{c1}) - {\check B}_L (u_{c2}) = {\check B}_L \left ( u_{c1} - u_{c2}, u_{c1} + u_{c2}\right ) 
\eeq

\noi so that by (\ref{3.43e}), by Lemma 3.3 and by (\ref{5.16e}) (\ref{5.17e})
\beq
\label{5.22e}
\parallel \omega^{\ell}  \partial_t \varphi_-(t) \parallel_2 \ \leq \ C\ b \ Z(1 - \ell n\ t)^{\beta}\  t^{-1- (1/2)(\ell + 1 - k)_+} \ h_2(t)
\eeq 

\noi for all $\ell > 0$. We can then impose the initial condition $\varphi_- (0) = 0$ and we obtain by integration in $(0, t]$
\beq
\label{5.23e}
\parallel \omega^{\ell}  \varphi_-(t) \parallel_2 \ \leq \ C\ b \ Z(1 - \ell n\ t)^{\beta}\  t^{- (1/2)(\ell + 1 - k)_+} \ h_2(t)\qquad \hbox{for $0 < \ell \leq 2$}
\eeq 

\noi by (\ref{5.14e}). We next define $\varphi_1 = \varphi_2 + \varphi_-$ so that $\varphi_1$ satisfies the equation 
\beq
\label{5.24e}
\partial_t \ \varphi_1 = {\check B}_L (u_{c1})\ .
\eeq

\noi Furthermore $\varphi_- = \varphi_1 - \varphi_2$ tends to zero when $t \to 0$ in the norms which appear in (\ref{5.23e}), and it follows from (\ref{5.20e}) (\ref{5.23e}) that 
\beq
\label{5.25e}
\parallel \omega^{\ell}  \varphi_i(t) \parallel_2 \ \leq \ C\ b^2 (1 - \ell n\ t)^{2\beta + 1} \qquad \hbox{for $0 < \ell \leq k$}
\eeq  

\noi for $i = 1, 2$, provided
\beq
\label{5.26e}
Z\ t^{-1/2}\ h_2(t) \leq b (1 - \ell n\ t)^{\beta + 1}
\eeq

\noi which can be ensured by taking $t$ sufficiently small because of  (\ref{5.15e}). We next define
\beq
\label{5.27e}
v_i = u_{ci} \exp \left ( i \varphi_i\right )
\eeq

\noi so that $(v_i, \varphi_i)$ satisfy the system (\ref{2.21e})  (\ref{2.22e}) since ${\check B}_L (u_{ci}) = {\check B}_L (v_{i})$. Furthermore, it follows from (\ref{5.16e}) (\ref{5.25e}) and Lemma 3.4, part (1), with $m = k$ that 
\bea
\label{5.28e}
\parallel v_{i} (t) ; V^k \parallel  &\leq& C\  b (1 - \ell n\ t)^{\beta} \left ( 1 + b^2(1 - \ell n\ t)^{2\beta + 1}\right )^2 \nn \\
&\leq& C\  b (1 + b^2)^2  (1 - \ell n\ t )^{5\beta + 2} \  .
\eea

\noi We next estimate $\parallel v_1 - v_2 \parallel_2$. From 
\beq
\label{5.29e}
v_1 - v_2 = \left ( u_{c1} - u_{c2}\right ) \exp \left ( i \varphi_1 \right ) + u_{c2} \left ( \exp \left ( i \varphi_-\right ) - 1 \right ) \exp \left ( i \varphi_2 \right ) \ ,
\eeq

\noi from (\ref{5.16e}) (\ref{5.17e}) (\ref{5.23e}) it follows that 
\beq
\label{5.30e}
\parallel v_{1} (t) - v_2(t)  \parallel_2\   \leq\  C\  Z(1 + b^2)  (1 - \ell n\ t)^{2\beta } \ h_2(t)\  .
\eeq

\noi Therefore $(v_i, \varphi_i)$ satisfy the assumptions of Proposition 5.1 with 
$$a = C\ b(1 + b^2)^2\quad , \quad \alpha = 5 \beta + 2\quad , \quad Y = C\ Z(1 + b^2)\quad , \quad h_1 = h_2 (1 - \ell n\ t)^{2\beta} $$

\noi so that the conditions (\ref{5.3e}) (\ref{5.4e}) reduce to (\ref{5.14e}) (\ref{5.15e}). It then follows from Proposition 5.1 that $(v_1, \varphi_1 ) = (v_2 , \varphi_2)$, so that $u_{c1} = u_{c2}$. \par \nobreak\hfill $\sq$\par

We now turn to the construction of solutions $(v, \varphi )$ of the system (\ref{2.21e})  (\ref{2.22e})  with prescribed asymptotic behaviour $(v_a, \varphi_a)$ as $t \to 0$, with $\varphi_a$ defined by (\ref{2.25e}) and $\varphi_a(1) = 0$, or equivalently of solutions of the system (\ref{2.30e})  (\ref{2.31e}) with $(w , \psi )$ tending to zero in a suitable sense when $t\to 0$. We first collect some preliminary estimates of $(G, \psi )$, of $(B_a, \varphi_a)$ and of $H_1$. \\

\noi {\bf Lemma 5.1.}\par

{\it (1) Let $v_a$, $w \in V^k$ and let $\psi$ satisfy (\ref{2.30e}). Then
\beq
\label{5.31e}
\parallel \omega^{\ell}  G \parallel_2 \ \leq \ C\parallel  w; H^k \parallel \ \parallel  2v_a + w ;H^k \parallel  \qquad \hbox{for $0 < \ell \leq k+ 1$}\ , 
\eeq  
\beq
\label{5.32e}
\parallel \omega^{\ell}  {\check G} \parallel_2 \ \leq \ C\parallel  w; V^k \parallel \ \parallel  2v_a + w ;V^k \parallel t^{-1}  \qquad \hbox{for $0 < \ell \leq 2k- 1$}\ , 
\eeq
\beq
\label{5.33e}
\parallel \omega^{\ell}  \partial_t \psi  \parallel_2 \ \leq \ C\parallel  w; V^k \parallel \ \parallel  2v_a + w ;V^k \parallel t^{-1-(\ell /2 + 1/2 - k)_+}  \qquad \hbox{for $\ell > 0$}\ .
\eeq 

(2) Let $v_a \in V^{k+1}$ and let  $\varphi_a$ satisfy (\ref{2.25e}). Then}  
\beq
\label{5.34e}
\parallel \omega^{\ell}  B_a  \parallel_2 \ \vee \ t \parallel \omega^{\ell} {\check B}_a \parallel_2\  \leq \ C\parallel  v_a; V^{k+1} \parallel^2  \qquad \hbox{\it for $0 < \ell \leq k+2$}\ , 
\eeq
\beq
\label{5.35e}
\parallel \omega^{\ell}  \partial_t \varphi_a \parallel_2 \ \leq \ C\parallel  v_a; V^{k+1} \parallel^2  \  t^{-1-(1/2)(\ell -2 - k)_+}  \qquad \hbox{\it for $\ell > 0$}\ .
\eeq 
\vskip 5 truemm

\noi {\bf Proof.} The estimates (\ref{5.31e}) (\ref{5.32e}) (\ref{5.34e}) follow from Lemma 3.3. The estimates (\ref{5.33e}) (\ref{5.35e}) follow from (\ref{3.43e}) and (\ref{5.32e}) (\ref{5.34e}).\par\nobreak\hfill $\sq$ \par

From now on we shall assume that $v_a$ satisfies the assumption
\beq
\label{5.36e}
v_a \in \left ( {\cal C} \cap L^{\infty}\right ) \left ( (0,1], V^{k+1}\right ) \ ,
\eeq

\noi we define
\beq
\label{5.37e}
a = \ \parallel v_a ;L^{\infty}\left ( (0, 1], V^{k+1}\right ) \parallel
\eeq

\noi and we assume that $\varphi_a$ satisfies (\ref{2.25e}) and $\varphi_a (1) = 0$, so that by integration of (\ref{5.35e})
\beq
\label{5.38e}
\parallel \omega^{\ell}  \varphi_a \parallel_2 \ \leq \ C\ a^2 \left ( t^{-(1/2)(\ell - 2 - k)_+} - \ell n\ t \right ) \qquad \hbox{for $\ell > 0$}
\eeq

\noi for all $t \in (0, 1]$.\par

We next estimate $H_1v_a$, where $H_1$ is defined by (\ref{2.32e}). We rewrite $H_1$ as
\beq
\label{5.39e}
H_1 = i L \cdot \nabla + M
\eeq

\noi where
\bea
\label{5.40e}
&&M = L \cdot K_a + (i/2) \nabla \cdot \sigma + (1/2) L^2 + {\check G}_S + g_{\not=}\ , \\
&&g_{\not=} = g \left ( w, 2v_a + w \right ) \ . \nn
\eea

\noi We shall use the auxiliary norm 
\beq
\label{5.41e}
\parallel f \parallel_{\star}\ = \ \parallel f\parallel_{\infty}\ \vee \ \parallel \nabla f \parallel_2 \ \vee\ \parallel \omega^k f \parallel_2
\eeq

\noi which satisfies
$$\parallel f_1 f_2  \parallel_{\star}\ \leq  \ \parallel f_1\parallel_{\star }\ \parallel f_2\parallel_{\infty}\ +  \ \parallel f_1 \parallel_{\infty} \ \parallel f_2  \parallel_{\star}$$

\noi by Lemma 3.1. \\

\noi {\bf Lemma 5.2.} {\it Let $v_a$ satisfy (\ref{5.36e}) with $a$ defined by (\ref{5.37e}). Let $I \subset (0, 1]$ and let $w \in ({\cal C} \cap L^{\infty}) (I, V^k)$, $\sigma \in {\cal C} (I, L^{\infty} \cap \dot{H}_1 \cap \dot{H}^{k+1} )$ with 
\beq
\label{5.42e}
\parallel w ; L^{\infty} (I, V^k )\parallel \ \leq C\ a \ .
\eeq

\noi Then the following estimate holds for all $t \in I$~:}
\bea
\label{5.43e}
&&\parallel H_1v_a; V^k \parallel\ \leq C\ a \Big  ( \parallel \sigma \parallel_{\star} \left ( 1 +\ \parallel \sigma \parallel_{\infty} \ + \ a^2 (1 - \ell n\ t)\right ) + \ \parallel \nabla \cdot \sigma \parallel_{\star}\nn \\
&&\qquad + \ a \parallel w; V^k \parallel \left ( t^{-1+(k-1)/2} + a^2 (1 - \ell n\ t )\right ) \Big ) \ .
\eea
\vskip 5 truemm

\noi {\bf Proof.} By Lemma 3.1, we estimate
$$\parallel H_1 v_a \parallel_2\ \leq \ \parallel L \parallel_{\infty}  \ \parallel \nabla v_a \parallel_2 \ + \ \parallel M \parallel_{\infty} \ \parallel v_a\parallel_2$$

$$\parallel xH_1 v_a \parallel_2\ \leq \ \parallel L \parallel_{\infty}  \ \parallel x\nabla v_a \parallel_2 \ + \ \parallel M \parallel_{\infty} \ \parallel xv_a\parallel_2$$

$$\parallel \omega^k H_1 v_a \parallel_2\ \leq \ \parallel  \omega^kL \parallel_2 \  \parallel \nabla v_a \parallel_{\infty} \ + \ \parallel L \parallel_{\infty} \ \parallel \omega^{k+1} v_a\parallel_2$$
$$+ \ \parallel \omega^k M \parallel_2\ \parallel  v_a \parallel_{\infty} \ + \ \parallel M \parallel_{\infty} \ \parallel \omega^k v_a\parallel_2$$

$$\parallel \omega^{k-1} x H_1 v_a \parallel_2\ \leq \left (  \parallel  L \parallel_{\infty}\ + \   \parallel \nabla L \parallel_2\right )  \parallel \omega^{k-1} \nabla v_a\parallel_2$$
$$+ \left (  \parallel M \parallel_{\infty}\ +\ \parallel \nabla M  \parallel_2 \right ) \parallel \omega^{k-1} v_a\parallel_2$$

\noi so that
\bea
\label{5.44e}
\parallel H_1 v_a;V^k \parallel &\leq& \left ( \parallel L\parallel_{\star}\ + \ \parallel M \parallel_{\star} \right ) \parallel v_a;V^{k+1} \parallel \nn\\
&\leq& a\left ( \parallel L\parallel_{\star}\ + \ \parallel M \parallel_{\star} \right ) \ .
\eea

\noi We next estimate
\beq
\label{5.45e}
\parallel L\parallel_{\star} \ \leq \ \parallel \sigma \parallel_{\star}\ + \ \parallel G\parallel_{\star} \ ,
\eeq
\beq
\label{5.46e}
\parallel M\parallel_{\star} \ \leq \ \parallel L\parallel_{\star}\left ( \parallel K_a\parallel_{\star} \ + \ \parallel L\parallel_{\infty}\right ) \ + \ \parallel \nabla \cdot \sigma \parallel_{\star}\ + \ \parallel {\check G}_S\parallel_{\star}\ + \ \parallel g_{\not=} \parallel_{\star} \ .
\eeq

\noi By Lemma 5.1, (\ref{3.42e}) and (\ref{5.37e}), we estimate
\beq
\label{5.47e}
\parallel G\parallel_{\star}  \ \leq \ C\ a \parallel w; V^k \parallel \ \leq C\ a^2 \ , 
\eeq
\beq
\label{5.48e}
\parallel {\check G}_S\parallel_{\star} \ \leq \ C\ a \parallel w; V^k \parallel \ t^{-1+(k-1)/2} \ , 
\eeq
\beq
\label{5.49e}
\parallel K_a\parallel_{\star}\ \leq \ C\ a^2(1 - \ell n\ t)\ .
\eeq

\noi Substituting (\ref{5.47e})-(\ref{5.49e}) into (\ref{5.45e}) (\ref{5.46e}) and substituting the result into (\ref{5.44e}) yields (\ref{5.43e}).\par \nobreak \hfill $\sq$ \par

We next give some estimates of solutions $w'$ of the linearized equation (\ref{2.34e}) associated with some $w \in X(I)$ where $X(I)$ is defined by (\ref{3.41e}) and $I = (0, \tau ]$ for some $\tau$, $0 < \tau \leq 1$. Such a $w$ satisfies
\beq
\label{5.50e}
\parallel  w(t) ; V^k \parallel \ \leq Y\ h(t)
\eeq

\noi for some $Y > 0$ and all $t \in I$. The following lemma is a variant of Lemma 4.1.\\

\noi {\bf Lemma 5.3.} {\it Let $v_a$ satisfy (\ref{5.36e}) with $a$ defined by (\ref{5.37e}) and let $\varphi_a$ be defined by (\ref{2.25e}) and $\varphi_a (1) = 0$. Let $0 < \tau \leq 1$ and $I = (0, \tau ]$. Let $w \in X(I)$ satisfy (\ref{5.50e}) for some $Y > 0$ and all $t \in I$, and let $\tau$ be sufficiently small so that
\beq
\label{5.51e}
Y\overline{h}(\tau ) \leq a \ .
\eeq

\noi Let $\psi$ be defined by  (\ref{2.30e}) with $\psi (0) = 0$. Then $\psi$ satisfies the estimate 
\beq
\label{5.52e}
\parallel \omega^{\ell} \psi \parallel_2\ \leq C\ a\ Y\ t^{-(\ell /2 + 1/2 - k)_+}\ h(t) \qquad \hbox{for $0 < \ell \leq k+2$}\ .
\eeq

\noi Let $w' \in {\cal C} (I, V^k)$ be a solution of the equation (\ref{2.34e}). Then the following estimates hold~:
\beq
\label{5.53e}
\Big | \partial_t \parallel w' \parallel_2 \Big  | \ \leq \ \parallel R_1 \parallel_2\ ,
\eeq
\beq
\label{5.54e}
\Big | \partial_t \parallel xw' \parallel_2 \Big  | \ \leq \ C\left ( 1 + a^2 (1 - \ell n\ t)\right ) \parallel w';V^k  \parallel\ + \ \parallel xR_1 \parallel_2\ ,
\eeq
\beq
\label{5.55e}
\left | \partial_t \parallel \omega^k w' \parallel_2 \right | \ \leq \ C\ a^2 ( 1 + a^2) (1 - \ell n\ t)^2  \parallel w';V^k  \parallel\ + \ \parallel \omega^k R_1 \parallel_2\ ,
\eeq
\beq
\label{5.56e}
\left | \partial_t \parallel \omega^{k-1} x w' \parallel_2 \right | \ \leq \ C\left ( 1 +  a^2 (1 - \ell n\ t)\right )^2   \parallel w';V^k  \parallel\ + \ \parallel \omega^{k-1} x  R_1 \parallel_2\ ,
\eeq

\noi where $R_1 = R - H_1 v_a$.}\\

\noi {\bf Proof.} It follows from (\ref{5.33e}) (\ref{5.50e}) (\ref{5.51e}) (\ref{3.39e}) that $\parallel \omega^{\ell}  \partial_t \psi  \parallel_2$ is integrable at $t = 0$ for $0 < \ell \leq k+2$, which gives a meaning to the assumption $\psi (0) = 0$. The estimate (\ref{5.52e})  then follows from (\ref{5.33e}) by integration. Furthermore (\ref{5.51e}) implies
\beq
\label{5.57e}
\parallel \omega^{\ell} \psi \parallel_2\ \leq C\ a\ Y\ \overline{h}(\tau ) \leq C\ a^2 \qquad \hbox{for $0 < \ell \leq k+2$}\ .
\eeq

\noi The proof of (\ref{5.53e})-(\ref{5.56e}) is a variant of that of the estimates (\ref{4.2e})-(\ref{4.4e}) of $v'$ in Lemma 4.1. We estimate in particular
\bea
\label{5.58e}
&&\parallel v; V^k \parallel\ \leq C\ a \  , \nn \\
&&\parallel \omega^{\ell} s \parallel_2\ \leq \ \parallel \omega^{\ell} s_a  \parallel_2\ + \ \parallel \omega^{\ell} \sigma \parallel_2\ \leq C\ a^2 (1 - \ell n\ t)
\eea

\noi for $0 \leq \ell \leq k+1$ by (\ref{5.38e}) (\ref{5.57e}). Furthermore
\beq
\label{5.59e}
\parallel \omega^{\ell} B \parallel\ \leq C\ a^2\qquad \hbox{for $0 < \ell \leq k+1$}\ ,
\eeq
\beq
\label{5.60e}
\parallel \omega^{\ell} g \parallel\ \leq C\ a^2\qquad \hbox{for $0 \leq \ell \leq k$}
\eeq

\noi by (\ref{5.36e}) (\ref{5.50e}) (\ref{5.51e}) and Lemma 3.3, while
$$\parallel \omega^{\ell} {\check B}_S \parallel\ \leq \ t^{(k+2 - \ell )/2} \parallel \omega^{k+2}  {\check B}_a \parallel_2 \ + \ t^{k-\ell /2 - 1/2} \parallel \omega^{2k-1} {\check G} \parallel_2$$
$$\leq C \left ( a^2\ t^{-1+(k+2-\ell )/2} + a\ Y\ t^{-1+k- \ell /2 - 1/2} \ h \right )$$

\noi for $0 < \ell \leq 2k-1$ by (\ref{3.42e}) (\ref{5.32e}) (\ref{5.34e}), so that by (\ref{5.51e})
\beq
\label{5.61e}
\parallel \omega^{\ell}  {\check B}_S  \parallel\ \leq C\ a^2\qquad \hbox{for $0 < \ell \leq k$}\ .
\eeq

\noi Substituting (\ref{5.58e})-(\ref{5.61e}) into the analogues for $w'$ of the estimates in the proof of Lemma 4.1, in particular (\ref{4.7e}) (\ref{4.8e}) yields (\ref{5.53e})-(\ref{5.56e}).\par\nobreak\hfill $\sq$ \par

We can now state the existence results of solutions of the linearized equation (\ref{2.34e}). \\

\noi {\bf Proposition 5.3.} {\it Let $1 < k < 2$. Let $v_a$ satisfy (\ref{5.36e}) with $a$ defined by (\ref{5.37e}) and let $\varphi_a$ be defined by (\ref{2.25e}) and $\varphi_a (1) = 0$. Let $0 < \tau \leq 1$ and $I = (0, \tau ]$. Let $w \in X(I)$ satisfy (\ref{5.50e}) for some $Y > 0$ and all $t \in I$ and let $\tau$ be sufficiently small to ensure (\ref{5.51e}). Let $R$ be defined by (\ref{2.33e}) and satisfy
\beq
\label{5.62e}
\parallel R; L^1 ((0, t], V^k ) \parallel\ \leq r\ h(t)
\eeq

\noi for some $r > 0$ and all $t \in I$. Then there exists a unique solution $w'\in X(I)$ of the equation (\ref{2.34e}) and $w'$ satisfies the estimate  
\beq
\label{5.63e}
\parallel w'(t) ; V^k \parallel\ \leq \left ( 1 + C(a) t (1 - \ell n\ t)^2 \right ) \left ( r + C(a) Y\ t^{(k-1)/2}\right ) h(t)
\eeq

\noi for some constant $C(a)$ depending on $a$ and for all $t \in I$.}\\

\noi {\bf Proof.} Let $0 < t_0 < \tau$ and let $w'_{t_0}$ be the solution of (\ref{2.34e}) in ${\cal C} (I, V^k)$ with initial condition $w'_{t_0}(t_0) = 0$. That solution is obtained by a minor variation of Proposition 4.1 including the inhomogeneous term $R_1$. We shall construct $w'$ as the limit of $w'_{t_0}$ when $t_0 \to 0$ and for that purpose we need estimates of $w'_{t_0}(t)$ in $V^k$ for $t \in [t_0, \tau ]$ that are uniform in $t_0$. From Lemma 5.2, especially (\ref{5.43e}) and from (\ref{5.50e}) (\ref{5.52e}), we obtain
\beq
\label{5.64e}
\parallel H_1 v_a ; V^k \parallel \ \leq h_1(t) \equiv C\ a^2 \ Y \left ( 1 + a^2 \ t^{1/2}(1 - \ell n\ t)\right ) t^{-(3-k)/2}\ h(t)\ .
\eeq

\noi On the other hand, from Lemma 5.3, especially (\ref{5.53e})-(\ref{5.56e}) and from (\ref{5.62e}) (\ref{5.64e}), we obtain 
\beq
\label{5.65e}
\parallel w'_{t_0}(t) ; V^k \parallel\ \equiv y(t) \leq \int_{t_0}^t f_1(t') \ y(t')\ dt' + f_2 (t)
\eeq

\noi where 
\beq
\label{5.66e}
f_1(t) = C \left ( 1 + a^2 (1 - \ell n\ t)\right )^2
\eeq
\bea
\label{5.67e}
&&f_2(t) = r\ h(t) + \int_0^t dt'\ h_1(t')\nn\\
&&\leq \left ( r + C\ a^2\ Y \left ( 1 + a^2\ t^{1/2} (1 - \ell n\ t) \right ) t^{(k-1)/2}\right ) h(t)
\eea

\noi by (\ref{3.39e}) (\ref{5.64e}). Integrating (\ref{5.65e}) with $y(t_0) = 0$ yields 
\beq
\label{5.68e}
y(t) \leq \int_{t_0}^t dt'\ f_1 (t') \ f_2(t') \exp \left ( \int_{t'}^t dt''\ f_1(t'')\right ) + f_2(t)\ .
\eeq

\noi Substituting (\ref{5.66e}) (\ref{5.67e}) into (\ref{5.68e}) yields
\bea
\label{5.69e}
&&y(t) \leq \left ( 1 + C \left ( 1 + a^2 (1 - \ell n\ t)\right )^2 t\right ) \exp \left ( C \left (1 + a^2 (1 - \ell n\ t)\right )^2 t \right ) \nn\\
&&\times \left ( r + C\ a^2\ Y \left ( 1 + a^2 \ t^{1/2} (1 - \ell n \ t)\right ) t^{(k-1)/2} \right ) h(t) \leq C(a, Y) h(t)
\eea

\noi uniformly in $t_0$. That estimate is of the form of (\ref{5.63e}). We now take the limit of $w'_{t_0}$ when $t_0 \to \infty$. Let $0 < t_0 \leq t_1 \leq \tau$. From the conservation of the $L^2$ norm of the difference of two solutions of (\ref{2.34e}), it follows that 
\beq
\label{5.70e}
\parallel w'_{t_0}(t) - w'_{t_1}(t)   \parallel_2 \ = \ \parallel w'_{t_0}(t_1)   \parallel_2 \ \leq C(a, Y) h(t_1)
\eeq

\noi for all $t \in [t_1, \tau ]$. It follows from (\ref{5.70e}) that $w'_{t_0}$ converges in $L_{loc}^{\infty}(I, L^2)$ to a limit $w' \in {\cal C} (I, L^2)$. From that convergence and from the uniform estimate (\ref{5.69e}), it follows that $w' \in {\cal C}(I, H^{k'}) \cap ({\cal C}_w \cap L^{\infty}) (I, V^k)$ for $0 \leq k' < k$ and that $w'_{t_0}$ converges to $w'$ in $L_{loc}^{\infty} (I, H^{k'})$ norm and weakly in $V^k$ pointwise in time. From the previous convergences and from the uniform estimate (\ref{5.69e}) of $w'_{t_0}$, it follows that $w'$ satisfies the same estimate in $I$. Clearly $w'$ is a solution of (\ref{2.34e}). It then follows from an inhomogeneous extension of Proposition 4.1 that $w' \in X (I)$. Finally the estimate (\ref{5.63e}) is a simplified version of (\ref{5.69e}). \par \nobreak \hfill $\sq$ \par

We can now derive the existence of solutions of the nonlinear system (\ref{2.30e}) (\ref{2.31e}). \\

\noi {\bf Proposition 5.4.} {\it Let $1 < k < 2$. Let $v_a$ satisfy (\ref{5.36e}) with $a$ defined by (\ref{5.37e}) and let $\varphi_a$ be defined by (\ref{2.25e}) and $\varphi_a (1) = 0$. Let $R$ be defined by (\ref{2.33e}) and satisfy (\ref{5.62e}) for all $t \in (0, 1]$. Then there exists $\tau$, $0 < \tau \leq 1$, depending on $(a,r)$ and there exists a unique solution $w \in X(I)$ of the equation (\ref{2.31e}) with $\psi$ satisfying (\ref{2.30e}) and $\psi (0) = 0$, where $I=(0, \tau ]$. In particular $w$ satisfies (\ref{5.50e}) for some $Y$ depending on $(a, r)$ and for all $t \in I$ and $\psi$ satisfies (\ref{5.52e}) for all $t \in I$.}\\

\noi {\bf Proof.} Let $0 < \tau \leq 1$. For $\tau$ sufficiently small, Proposition 5.3 defines a map $\Gamma : w \to w'$ from $X(I)$ into itself. We shall show that for $\tau$ sufficiently small, the map $\Gamma$ is a contraction on the subset ${\cal R}$ of $X(I)$ defined by (\ref{5.50e}) for a suitable choice of $Y$ in the norm considered in Proposition 5.1. \par

We first ensure that ${\cal R}$ is stable under $\Gamma$. From (\ref{5.63e}) it follows that
\beq
\label{5.71e}
\parallel w';X(I)\parallel\ \leq \left ( 1 + C(a)\ \tau (1 - \ell n\ \tau )^2\right ) \left ( r + C(a) Y \tau^{(k-1)/2}\right )
\eeq

\noi and this can be made smaller than $Y$ by taking $Y = 2r$ and $\tau$ sufficiently small. We next show that $\Gamma$ is a contraction on ${\cal R}$ for the $L^2$ norm of $w$. Let $w_i \in {\cal R}$ and $w'_i = \Gamma w_i$, $i = 1,2$, let $w_{\pm} = (1/2) (w_1 \pm w_2)$ and similarly for $w'_i$. All those quantities belong to ${\cal R}$. We define the norms
\bea
\label{5.72e}
&&Y_- = \ \mathrel{\mathop {\rm Sup}_{t\in I}}\ h(t)^{-1} \parallel w_-(t) \parallel_2 \\
&&Y'_- = \ \mathrel{\mathop {\rm Sup}_{t\in I}}\ h(t)^{-1} \parallel w'_-(t) \parallel_2 
\label{5.73e}
\eea

\noi and we estimate $Y'_-$ in terms of $Y_-$ by Lemma 4.2 with $v_- = w_-$, $v_+ = v_a + w_+$, $v'_+ = v_a + w'_+$, $s_- = \sigma_-$, $s_+ = s_a + \sigma_+$ and $y$ defined by (\ref{4.13e}). From (\ref{5.37e}) (\ref{5.50e}) (\ref{5.51e}) it follows that
\beq
\label{5.74e}
y \ \vee\ \parallel v_+; V^k \parallel \ \vee \ \parallel v'_+; V^k \parallel \ \leq C\ a \ .
\eeq

\noi From (\ref{5.38e}) (\ref{5.52e}) (\ref{5.51e}) it follows that  
\beq
\label{5.75e}
\parallel  \omega^{\ell} s_+ \parallel _2 \ \leq C\ a^2(1 - \ell n\ t) \qquad \hbox{for $0 \leq \ell \leq k+1$}\ .
\eeq

\noi From (\ref{4.15e}) it follows that
\beq
\label{5.76e}
\parallel  \omega^{\ell} \sigma_- \parallel _2 \ \leq C\ a \ Y \ t^{-(\ell +2 - k)/2}\ h(t)  \qquad \hbox{for $0 \leq \ell \leq 1$}\ .
\eeq

\noi Substituting (\ref{5.74e})-(\ref{5.76e}) into (\ref{4.14e}) yields
\beq
\label{5.77e}
\left | \partial_t \parallel w'_-  \parallel_2 \right | \ \leq \ C\ a^2 \ Y_- \left ( 1 +  a^2 \ t^{1/2} (1 - \ell n\ t)\right ) t^{-(3-k)/2}\ h(t)
\eeq

\noi and therefore by integration over time
\beq
\label{5.78e}
Y'_- \leq C\ a^2\ Y_- \left ( 1 + a^2 \tau^{1/2} (1 - \ell n\ \tau )\right ) \tau^{(k-1)/2}
\eeq

\noi which implies the contraction property for the norm (\ref{5.72e}) (\ref{5.73e}) for $\tau$ sufficiently small. The existence result now follows from the fact that ${\cal R}$ is closed for that norm.\par

Finally the uniqueness result follows from Proposition 5.1 with $\alpha = 0$. \par \nobreak\hfill $\sq$ \par

So far we have solved the auxiliary system (\ref{2.30e}) (\ref{2.31e}) for general $v_a$ satisfying (\ref{5.36e}) and under the decay assumption (\ref{5.62e}) on the remainder $R$ defined by (\ref{2.33e}). We now take $v_a = U(t) v_+$. By (\ref{3.8e}), that $v_a$ satisfies (\ref{5.36e}) for $v_+ \in V^{k+1}$ with 
\beq
\label{5.79e}
a = \ \parallel  v_a; L^{\infty} ((0, 1], V^{k+1})\parallel \ \leq \ 2\parallel v_+ ; V^{k+1} \parallel \ .
\eeq

\noi We now show that the corresponding remainder satisfies (\ref{5.62e}) with $h(t) =$\break\noindent $t(1 - \ell n\ t)^2$.\\

\noi {\bf Proposition 5.5.} {\it Let $1 < k < 2$. Let $v_+ \in V^{k+1}$ and $v_a (t) = U(t) v_+$. Let $R$ and $\varphi_a$ be defined by (\ref{2.33e}) (\ref{2.25e}) and $\varphi_a (1) = 0$. Then $R \in {\cal C} ((0, 1], V^k)$ and $R$ satisfies the estimate
\beq
\label{5.80e}
\parallel  R(t) ; V^k \parallel \ \leq C\ a^3 (1 + a^2) (1 - \ell n\ t)^2 
\eeq

\noi for all $t \in (0, 1]$.}\\

\noi {\bf Proof.} For $v_a = U(t) v_+$, the remainder $R$ takes the form 
\beq
\label{5.81e}
R = - i K_a \cdot \nabla v_a - \left ( (i/2) \nabla \cdot s_a + (1/2) K_a^2 + {\check B}_{aS} + g_a \right ) v_a \ .
\eeq

\noi Using the fact that the norm $\parallel\ \cdot \ \parallel_{\star}$ defined by (\ref{5.41e}) satisfies
\beq
\label{5.82e}
\parallel f_1 f_2 ; V^k \parallel \ \leq \ C \parallel f_1 \parallel_{\star} \ \parallel f_2; V^k \parallel
\eeq

\noi by Lemma 3.1, we estimate
$$\parallel R; V^k \parallel\ \leq \ C \Big ( \parallel K_a \parallel_{\star} \ \parallel v_a ; V^{k+1} \parallel\ + \Big ( \parallel \nabla \cdot s_a \parallel_{\star} \ + \ \parallel K_a \parallel_{\star}^2$$
\beq
\label{5.83e}
+\ \parallel {\check B}_{aS} \parallel_{\star}\ + \ \parallel g_a \parallel_{\star}\Big ) \parallel v_a ; V^k \parallel \ .
\eeq

\noi We then estimate
\bea
\label{5.84e}
&&\parallel B_a \parallel_{\star} \ \leq C\ a^2 \\
&&\parallel {\check B}_{aS} \parallel_{\star}\ \leq \left ( 1 + C\ t^{(k-1)/2}\right ) \parallel \omega^{k+2}{\check B}_a \parallel_2\ \leq C\ a^2
\label{5.85e}
\eea

\noi by (\ref{5.34e}) (\ref{3.42e}),
\beq
\label{5.86e}
\parallel  s_a \parallel_{\star} \ \vee \ \parallel \nabla \cdot s_a \parallel_{\star}\ \leq C\ a^2 (1 - \ell n\ t)
\eeq

\noi by (\ref{5.38e}) and 
\beq
\label{5.87e}
\parallel g_a \parallel_{\star}\ \leq C\ a^2 
\eeq

\noi by Lemma 3.3. Substituting (\ref{5.84e})-(\ref{5.87e}) into (\ref{5.83e}) yields (\ref{5.80e}). The continuity of $R$ in $V^k$ follows immediately from that of $v_a$ in $V^{k+1}$ and from the estimates.\par\nobreak\hfill $\sq$ \par

Proposition 5.5 implies that $R$ satisfies the assumption (\ref{5.62e}) with 
\beq
\label{5.88e}
h(t) = t (1 - \ell n\ t)^2
\eeq

\noi so that 
\beq
\label{5.89e}
\overline{h} (t) = t^{(k-1)/2} (1 - \ell n\ t)^2\ .
\eeq

\noi Putting together Propositions 5.4 and 5.5, we obtain the final result for the Cauchy problem at time zero for the system (\ref{2.21e}) (\ref{2.22e}) in the following form.\\

\noi {\bf Proposition 5.6.} {\it Let $1 < k < 2$. Let $v_+ \in V^{k+1}$ with 
\beq
\label{5.90e}
a_+ = \ \parallel  v_+ ; V^{k+1} \parallel \ .
\eeq

\noi Let $v_a = U(t) v_+$ and let $\varphi_a$ be defined by (\ref{2.25e}) and $\varphi_a (1) = 0$. Then there exists $\tau$, $0 < \tau \leq 1$, depending on $a_+$ and there exists a unique solution $(v, \varphi )$ of the system (\ref{2.21e}) (\ref{2.22e}) such that $v \in ({\cal C} \cap L^{\infty})(I, V^k)$, $\varphi \in {\cal C} (I, H_>^{k+2})$, $v-v_a \in X(I)$ and $(\varphi - \varphi_a)(0) = 0$, where $I = (0, \tau ]$ and $X(I)$ is defined by (\ref{3.41e}) with $h$ given by (\ref{5.88e}). The solution $(v , \varphi )$ satisfies the estimates
\beq
\label{5.91e}
\parallel v(t) - v_a (t) ; V^k \parallel\ \leq Y\ t(1 - \ell n\ t)^2
\eeq
\beq
\label{5.92e}
\parallel \omega^{\ell} \left ( \varphi (t) - \varphi_a (t)\right ) \parallel_2\ \leq C\ a_+\ Y\ t^{1 - (\ell /2 + 1/2 - k)_+} (1 - \ell n \ t)^2 \quad \hbox{\it for $0 < \ell \leq k+2$}
\eeq

\noi for some $Y > 0$ depending on $a_+$ and for all $t \in I$.}\\

We can now state the final result on the Cauchy problem for the system (\ref{1.1e}) with prescribed asymptotic behaviour at infinity.\\

\noi {\bf Proposition 5.7.} {\it Let $1 < k < 2$. Let $u_+ \in FV^{k+1}$ with 
\beq
\label{5.93e}
a_+ = \ \parallel  u_+ ; FV^{k+1} \parallel \ .
\eeq

\noi Let $u_a$ be defined by (\ref{2.35e}) (\ref{2.36e}) with $v_a = U(t) v_+$, $v_+ = \overline{Fu_+}$ and $\varphi_a$ defined by (\ref{2.25e}) and $\varphi_a (1) = 0$. Then there exists $T \geq 1$ depending on $a_+$ and there exists a unique solution $u$ of the system (\ref{1.1e}) such that $\widetilde{u} \in {\cal C} (I, FV^k)$, where $I = [T, \infty )$ and $\widetilde{u}  (t) = U(-t) u(t)$, and such that 
\beq
\label{5.94e}
\parallel  \widetilde{u} (t) - \widetilde{u} _a(t); FV^k \parallel \ \leq C(a_+) t^{-1} (1 + \ell n\ t)^{2+k}
\eeq

\noi for all $t \in I$. More precisely, $\widetilde{u}$ satisfies the estimates}
\beq
\label{5.95e}
\parallel |x|^{\ell} \widetilde{u}(t) \parallel _2 \ \leq C\ a_+ \left ( 1 + a_+^2 (1 + \ell n\ t)\right )^{\ell} \qquad \hbox{\it for $0 \leq \ell \leq k$}\ ,
\eeq
\beq
\label{5.96e}
\parallel |x|^{\ell} \nabla \widetilde{u}(t) \parallel _2 \ \leq C\ a_+ \left ( 1 + a_+^2 (1 + \ell n\ t)\right )^{\ell} \qquad \hbox{\it for $0 \leq \ell \leq k-1$}\ ,
\eeq
\beq
\label{5.97e}
\parallel  \widetilde{u} (t) - \widetilde{u} _a(t) \parallel_2\ \vee\ \parallel \nabla \left ( \widetilde{u}(t) - \widetilde{u} _a(t)\right ) \parallel _2 \ \leq C(a_+)t^{-1} (1 + \ell n\ t)^2\ , 
\eeq
\beq
\label{5.98e}
\parallel |x|^{k} \left ( \widetilde{u}(t) - \widetilde{u} _a(t)\right )\parallel _2 \ \leq C(a_+)  t^{-1} (1 + \ell n\ t)^{k+2}\ , 
\eeq
\beq
\label{5.99e}
\parallel |x|^{k-1} \nabla \left ( \widetilde{u}(t) - \widetilde{u} _a(t)\right )\parallel _2 \ \leq C(a_+)  t^{-1} (1 + \ell n\ t)^{3}\ . 
\eeq

\vskip 5 truemm

\noi {\bf Proof.} We first prove the existence of $u$ with the properties stated. Let $(v, \varphi )$ be the solution of the system (\ref{2.21e}) (\ref{2.22e}) obtained in Proposition 5.6 and define $u$ by (\ref{2.6e}) (\ref{2.16e}). Then $u$ is a solution of the system (\ref{1.1e}) defined in $I = [T , \infty )$ with $T = \tau^{-1}$. The properties of $u$ follow from those of $(v, \varphi )$ and from the estimates which we now derive. By (\ref{2.6e}) (\ref{2.35e}) and (\ref{3.8e}) it is sufficient to estimate $u_c$ and $u_c - u_{ca}$ in $V^k$. From (\ref{5.38e}) and from Lemma 3.4, part (1) with $m = 2$ and part (3), it follows that 
\beq
\label{5.100e}
\parallel \omega^{\ell} u_{ca} (t) \parallel _2\ \leq C\ a_+ \left ( 1 + a_+^2 (1 - \ell n\ t)\right )^{\ell} \qquad \hbox{for $0 \leq \ell \leq k$}\ , 
\eeq
\beq
\label{5.101e}
\parallel \omega^{\ell} x u_{ca} (t) \parallel _2\ \leq C\ a_+ \left ( 1 + a_+^2 (1 - \ell n\ t)\right )^{\ell} \qquad \hbox{for $0 \leq \ell \leq k-1$}\ .
\eeq

\noi We next estimate the difference
\bea
\label{5.102e}
u_c - u_{ca} &=& v \exp (- i \varphi ) - v_a \exp (- i \varphi_a)\nn \\
&=&\left ( v (\exp (-i \psi ) - 1 ) + v - v_a\right ) \exp (- i \varphi_a)
\eea

\noi with $\psi = \varphi - \varphi_a$. From (\ref{5.91e}) (\ref{5.92e}) and Lemma 3.4, part (2), it follows that 
$$\parallel v(\exp (-i \psi ) - 1 ) ; V^k\parallel \ \leq C\ a_+ \left ( a_+\ Y\ t(1 - \ell n\ t)^2 + \left ( a_+ \ Y\ t(1 - \ell n \ t)^2\right )^2 \right )$$
\beq
\label{5.103e}
\leq C \ a_+^2 \ Y\ t (1 - \ell n\ t)^2
\eeq

\noi for $\tau$ sufficiently small so that $a_+ Y \tau (1 - \ell n \tau )^2 \leq 1$. From (\ref{5.91e}), (\ref{5.103e}) (\ref{5.38e}) and Lemma 3.4, part (1) with $m = 2$, it then follows that 
\beq
\label{5.104e}
\parallel <x> \left ( u_c (t) - u_{ca} (t) \right ) \parallel_2 \ \leq C \left ( 1 + a_+^2 \right ) Y\ t(1 - \ell n\ t)^2\ , 
\eeq
\beq
\label{5.105e}
\parallel \omega^k \left ( u_c (t) - u_{ca} (t) \right ) \parallel_2 \ \leq C \left ( 1 + a_+^2 \right )^{k+1}\  Y\ t(1 - \ell n\ t)^{k+2}\ , 
\eeq
\beq
\label{5.106e}
\parallel \omega^{k-1}x \left ( u_c (t) - u_{ca} (t) \right ) \parallel_2 \ \leq C \left ( 1 + a_+^2 \right )^2\  Y\ t(1 - \ell n\ t)^3\ .
\eeq

\noi The estimates (\ref{5.97e})-(\ref{5.99e}) follow from (\ref{5.104e})-(\ref{5.106e}) and the estimates (\ref{5.95e}) (\ref{5.96e}) follow from (\ref{5.100e}) (\ref{5.101e}) and (\ref{5.97e})-(\ref{5.99e}).\par

Uniqueness of $u$ follows from Proposition 5.2 with $h_2(t) = t(1 - \ell n\ t)^2$ and $\beta = k$, which satisfy the conditions (\ref{5.14e}) (\ref{5.15e}).\par \nobreak \hfill $\sq$\par

\section*{Acknowledgments}
\hspace*{\parindent} This work arose from discussions of one of us (G.V.) with Professor Jalal Shatah who suggested that the methods previously used for the Maxwell-Schr\"odinger system could be applied to the Modified Schr\"odinger Map. We are greatly indebted to Professor Shatah for that fruitful suggestion and for introducing us to the Modified Schr\"odinger Map. Part of this work was done while one of us (G. V.) was visiting the Courant Institute of Mathematical Sciences and the Institut des Hautes Etudes Scientifiques. He is very grateful to Professor Jalal Shatah and to Professor Jean-Pierre Bourguignon for the hospitality extended to him at those two institutions. We are grateful to Professor Patrick G\'erard for an enlightening conversation and to Professor Tohru Ozawa and Professor Yoshio Tsutsumi for correspondence.

\begin{appendix}
\mysection{- Appendix. Properties of $\Sigma^k$.}
\hspace*{\parindent} We work in arbitrary space dimension $n$, although we need only the special case $n=2$. We first show that $\Sigma^k \subset V^k$. For that purpose it suffices to prove the following lemma.\\

\noi {\bf Lemma A.1.} {\it Let $1 < k < 2$. Then the following estimate holds
\beq
\label{A.1e}
\parallel < \omega >^{k-1} xv \parallel_2 \ \leq C\left ( \parallel <\omega >^k v \parallel_2 \ + \ \parallel <x>^k v \parallel_2 \right )
\eeq

\noi or equivalently}
\beq
\label{A.2e}
\parallel \omega^{k-1} xv \parallel_2 \ \leq \ C\parallel \omega^k v \parallel_2^{1-1/k} \ \parallel |x|^k v \parallel_2^{1/k}\ . 
\eeq
\vskip 5 truemm

\noi {\bf Proof.} From the elementary estimate
\beq
\label{A.3e}
\parallel v \parallel_2 \ \leq C\left ( \parallel \omega^k v \parallel_2 \ + \ \parallel |x|^k v \parallel_2 \right )
\eeq

\noi it follows that (\ref{A.1e}) is equivalent to
\beq
\label{A.4e}
\parallel \omega^{k-1} xv \parallel_2 \ \leq \ C\left ( \parallel \omega^k v \parallel_2 \ +\ \parallel |x|^k v \parallel_2\right ) \ . 
\eeq

\noi Clearly (\ref{A.2e}) implies (\ref{A.4e}). Conversely (\ref{A.2e}) follows from (\ref{A.4e}) by a dilation of $v$ followed by an optimization of the dilation parameter.\par

In order to prove (\ref{A.1e}), we use a dyadic decomposition. Let $\widehat{\psi} \in {\cal C}^{\infty} ({I\hskip-1truemm R}^n, {I\hskip-1truemm R})$, $0 \leq \widehat{\psi} \leq 1$, $\widehat{\psi}(\xi ) = 1$ (resp. $0$) for $|\xi | \leq 1$ (resp. $\geq 2$), let $\widehat{\varphi}_0 = \widehat{\psi}$ and $\widehat{\varphi}_j (\xi ) = \widehat{\psi} (2^{-j} \xi ) - \widehat{\psi} (2^{-(j-1)}\xi ) = \widehat{\varphi}_1 (2^{-(j-1)} \xi )$,

\noi so that
$$\varphi_j (x) = 2^{(j-1)n} \ \varphi_1\left  (2^{j-1} \ x\right )$$

\noi and therefore
\beq
\label{A.5e}
\parallel |x|^{\ell} \ \varphi_j\parallel_1 \ = \ 2^{-(j-1)\ell}\  \parallel |x|^{\ell}\ \varphi_1 \parallel_1 \qquad \hbox{for $\ell \geq 0$}
\eeq

\noi where we use the notation $\widehat{\varphi} = F \varphi$. Clearly for $\ell \geq 0$
\beq
\label{A.6e}
C^{-1} \parallel < \omega >^{\ell } v \parallel_2^2 \ \leq \ \sum_{j\geq 0}\ 2^{2j\ell}  \parallel \varphi_j \star v \parallel_2^2 \ \leq  \ C\parallel <\omega >^{\ell} v \parallel_2^2 \ .
\eeq

\noi Let now $\ell = k-1$ so that $0 < \ell < 1$. We estimate
\beq
\label{A.7e}
\parallel < \omega >^{\ell } xv \parallel_2^2 \ \leq \ C \sum_{j\geq 0}\ 2^{2j\ell}  \parallel \varphi_j \star xv \parallel_2^2 \ .
\eeq

\noi Now
$$\varphi_j \star xv = x \left ( \varphi_j \star v\right ) - \left ( x \varphi_j\right ) \star v$$

\noi so that by (\ref{A.5e}) and the Young inequality
\beq
\label{A.8e}
\parallel \varphi_j \star xv \parallel_2\ \leq \ \parallel x \left ( \varphi_j \star v\right )  \parallel_2\ +\ 2^{-(j-1)} \parallel |x| \ \varphi_1 \parallel_1\ \parallel v\parallel_2 \ .
\eeq

\noi Substituting (\ref{A.8e}) into (\ref{A.7e}) yields
\beq
\label{A.9e}
\parallel < \omega >^{\ell } xv \parallel_2^2 \ \leq \ C \sum_{j\geq 0}\ 2^{2j\ell}  \parallel x\left ( \varphi_j \star v\right )  \parallel_2^2 \  + \ C \parallel v\parallel_2^2 \ .
\eeq

\noi We estimate the sum in the RHS by the H\"older inequality in $(j, x)$ as
$$\sum_{j\geq 0} 2^{2j\ell} \parallel x\left ( \varphi_j \star v\right )  \parallel_2^2\ \leq \left ( \sum_{j\geq 0} 2^{jk}  \parallel \varphi_j \star v  \parallel_2^2\right )^{\ell /k} \left ( \sum_{j\geq 0} \parallel |x|^k \left ( \varphi_j \star v\right )  \parallel_2^2\right )^{1 /k}$$ 
\beq
\label{A.10e}
\leq C\parallel  <\omega >^k v \parallel_2^{2\ell /k} \left ( \sum_{j\geq 0} \parallel |x|^k \left ( \varphi_j \star v\right )  \parallel_2^2\right )^{1 /k}
\eeq

\noi by (\ref{A.6e}). Now from the inequality
$$\left |  | x|^k - |y|^k\right | \leq k |x-y|\left ( |x-y|^{\ell} + |y|^{\ell} \right )$$

\noi we obtain
$$\left | | x|^k \left ( \varphi_j \star v\right ) - \varphi_j \star |x|^k v\right | \leq k\left ( \left (  |x|^k |\varphi_j \right ) \star |v| + \left ( |x|\ |\varphi_j|\right ) \star |x|^{\ell}  |v| \right )$$

\noi and therefore by (\ref{A.5e}) and the Young inequality
$$\parallel |x|^k \left ( \varphi_j \star v\right )\parallel_2\  \leq \ \parallel  \varphi_j \star |x|^k v\parallel \ + \  k\Big  (2^{-(j-1)k} \parallel |x|^k \varphi_1 \parallel_1 \ \parallel v \parallel_2$$
 \beq
\label{A.11e}
+ \ 2^{-(j-1)} \parallel |x| \varphi_1 \parallel_1 \ \parallel |x|^{\ell} v \parallel_2 \Big ) 
\eeq

\noi so that by (\ref{A.6e})
\beq
\label{A.12e}
\sum_{j\geq 0} \parallel |x|^k \left ( \varphi_j \star v\right )\parallel_2^2\  \leq \ C\left ( \parallel |x|^k v \parallel_2^2 \ +\ \parallel |x|^{\ell} v \parallel_2^2 \ + \ \parallel v \parallel_2^2 \right )\ .
\eeq

\noi Substituting (\ref{A.12e}) into (\ref{A.9e}) (\ref{A.10e}) yields (\ref{A.1e}). \par\nobreak\hfill $\sq$\par

We next prove that $\Sigma^k$ is stable under the free Schr\"odinger evolution $U(\ \cdot \ )$ and for that purpose we prove the estimate (\ref{3.11e}). We treat only the case $1 < k < 2$, but the extension to general $k$ is straightforward.\\

\noi {\bf Lemma A.2.} {\it Let $1 < k < 2$. Then} 
\beq
\label{A.13e}
\parallel U(t) v;\Sigma^k \parallel\ \leq \ C \left ( \parallel v; \Sigma^k \parallel\ + \ |t|^k \parallel \omega^k v \parallel_2\right )\ .
\eeq
\vskip 5 truemm

\noi {\bf Proof.} By (\ref{A.3e}), one can use for $\Sigma^k$ the equivalent norm 
$$\parallel v;\Sigma^k \parallel \ = \ \parallel \omega^k v \parallel _2\ +\ \parallel |x|^kv\parallel_2$$

\noi and by the commutation relation (\ref{3.6e}), it is sufficient to prove that 
\beq
\label{A.14e}
\parallel |x + it \nabla |^k v \parallel\ \leq \ C \left ( \parallel |x|^k v \parallel_2\ + \ |t|^k \parallel \omega^k v \parallel_2\right )\ .
\eeq

\noi By homogeneity (namely dilation of $v$ by $|t|^{1/2}$), (\ref{A.14e}) is equivalent to the special case $t= 1$. It is then sufficient to prove that 
\beq
\label{A.15e}
\parallel \omega^k \exp (i x^2/2) v \parallel_2 \ \leq \ C \left ( \parallel \omega^k v \parallel_2\ + \  \parallel |x|^k v \parallel_2\right )\ .
\eeq

\noi We first prove (\ref{A.15e}) with $k$ replaced by $\ell$ with $0 < \ell < 1$. By (\ref{3.16e})
\beq
\label{A.16e}
\parallel \omega^{\ell} \exp (i x^2/2) v \parallel_2^2 \ = \ C \int dy \ |y|^{-n-2\ell} \parallel \left ( \tau_y - \tau_{-y}\right ) \exp (i x^2/2) v \parallel_2^2 \ .
\eeq

\noi Now
$$\left ( \tau_y - \tau_{-y}\right ) \exp (i x^2/2) v = \left ( \tau_y  \exp (i x^2/2) \right ) \left ( \tau_y v-v\right )$$
$$- \ \tau_{-y} \exp (i x^2/2) \left ( \tau_{-y} v-v\right ) + 2i \exp \left ( (i (x^2 + y^2)/2\right ) \sin (xy) v(x)$$

\noi so that
\bea
\label{A.17e}
\parallel \left ( \tau_y - \tau_{-y}\right ) \exp (i x^2/2) v\parallel_2 &\leq& \parallel \tau_y v-v\parallel_2 \ +\ \parallel \tau_{-y} v-v \parallel_2\nn \\
&&+ \ 2\parallel \sin (\cdot \ y ) v ( \cdot ) \parallel_2
\eea

\noi Substituting (\ref{A.17e}) into (\ref{A.16e}) and using again (\ref{3.16e}) yields
$$\parallel \omega^{\ell} \exp (i x^2/2) v \parallel_2^2 \ \leq \ C \left ( \parallel \omega^{\ell} v\parallel_2^2\ + \ \int dy\  |y|^{-n-2\ell}\ \parallel \sin ( \cdot \ y) v(\cdot ) \parallel_2^2 \right ) \ .$$

\noi The last integral is
$$\int dy\ dx \ |y|^{-n-2\ell} \sin^2 (xy) |v(x)|^2 = \ C \parallel |x|^{\ell} v \parallel_2^2$$

\noi which completes the proof of (\ref{A.15e}) with $k$ replaced by $\ell$.\par

We next take $\ell = k-1$ and we estimate
\begin{eqnarray*}
&&\parallel |x + i \nabla |^k v\parallel_2\ = \ \parallel |x + i \nabla |^{\ell} (x + i \nabla ) v\parallel_2\\
&&\leq\ C\left (  \parallel |x|^k v \parallel_2 \ + \  \parallel \omega^{\ell} xv\parallel_2 \ + \ \parallel |x|^{\ell} \nabla v \parallel_2 \ + \ \parallel \omega^{k} v\parallel_2 \right ) \\
&&\leq\ C \parallel v; \Sigma^k \parallel
\end{eqnarray*}

\noi by (\ref{A.15e}) with $k$ replaced by $\ell$ and Lemma A.1.  \par\nobreak\hfill $\sq$\par

\end{appendix}

\end{document}